\documentclass[reqno,12pt]{article}
\usepackage{amsmath,amsthm,amsfonts,amscd}
\usepackage{xy}
\usepackage{graphicx}

\newtheorem{thm}{Theorem}[section]
\newtheorem{theorem}{Theorem}[subsection]
\newtheorem{lem}[thm]{Lemma}
\newtheorem{cor}[thm]{Corollary}
\newtheorem{prop}[theorem]{Proposition}
\newtheorem{lemma}[theorem]{Lemma}

\theoremstyle{definition}
\newtheorem{rem}[thm]{Remark}
\newtheorem{defn}[thm]{Definition}

\newcommand{\W}{\mathfrak{W}}

\begin{document}
\title{On the equivariant Gromov-Witten Theory of
$\mathbb{P}^{2}$-bundles over curves}
\author{Amin Gholampour}\maketitle
\begin{abstract}
We compute section class relative equivariant Gromov-Witten
invariants of the total space of $\mathbb{P}^{2}$-bundles of the
form
$$\mathbb{P}(\mathcal{O}\oplus L_{1}\oplus L_{2})\rightarrow C,$$
where $C$ is a genus $g$ curve, and $\mathcal{O}$ is the trivial
bundle, and $L_{1}$ (resp. $L_{2}$) is an arbitrary line bundle of
degree $k_{1}$ (resp. $k_{2}$) over $C$.

We prove a gluing formula for the partition functions of these
invariants. The gluing formula allows us to compute the partition
function in general case in terms of the basic partition functions
for the case of $g=0$, relative to one, two or three fibers. We
compute these basic partition functions via localization
techniques combined with relations arising from the gluing
formula. These give rise to explicit $3\times 3$  matrices $G$,
$U_{1}$ and $U_{2}$ with entries in
$\mathbb{Q}((u))(t_{0},t_{1},t_{2})$, where $u$ is the genus
parameter, and $t_{0},t_{1},t_{2}$ are the equivariant parameters.
Then we prove that the partition function of the section class,
ordinary equivariant Gromov-Witten invariants of $X$ is given by
(Theorem \ref{thm:matrix formula}):
$$\operatorname{tr}\left(G^{g-1}U_{1}^{k_{1}}U_{2}^{k_{2}}\right).$$

As an application, we establish a formula for the partition
function of the ordinary Gromov-Witten invariants of any
$\mathbb{P}^{2}$-bundle $X$ over a curve of genus $g$ for any
class which is a Calabi-Yau section class, i.e. a curve class
$\beta_{cs}$ such that $K_{X}\cdot \beta_{cs}=0$ and $F\cdot
\beta_{cs}=1$, where $K_{X}$ is the canonical bundle of $X$ and
$F$ is the class of the fiber. We prove that this partition
function is given by (Theorem \ref{thm:Calabi-Yau class}):
$$3^{g}\left(2\operatorname{sin}\frac{u}{2}\right)^{2g-2}.$$
\end{abstract}

\tableofcontents

\section{Introduction} \label{sec:notation}
Let $X$ be a $\mathbb{P}^{2}$-bundle over a curve $C$ of genus
$g$, and we denote the cohomology class of the fiber by $F$.
\begin{defn} \label{defn:section class}
A class $\beta \in H^{4}(X,\mathbb{Z})$ is called \textbf{section
class} if $$F \cdot \beta=1.$$

We say $\beta$ is a \textbf{Calabi-Yau class} if
$$K_{X}\cdot \beta=0.$$ $\beta$ is called \textbf{Calabi-Yau
section class} if both conditions above hold.
\end{defn}

\begin{rem} \label{rem:general section class}
A section class is not necessarily represented by a geometric
section of the bundle $X$. It could be for example a section with
a number of fiber curves (curves which are included in the fibers
of $X$) are attached to it.
\end{rem}

If $\beta$ is a Calabi-Yau class then the virtual dimension of the
moduli space of degree $\beta$, genus $h$ stable maps to $X$ is
zero:
\[\text{virdim }
\overline{M}_{h}(X,\beta)=0. \]

Let $\beta_{cs} \in H^{4}(X,\mathbb{Z})$ be a Calabi-Yau section
class. The \textbf{partition function} of the class $\beta_{cs}$
Gromov-Witten invariants of $X$ by
\[
Z_{\beta_{cs}}(g)=\sum_{h=0}^{\infty}u^{2h-2}
\int_{[\overline{M}_{h}(X,\beta_{cs})]^{vir}}1,
\]
where $$[\overline{M}_{h}(X,\beta_{cs})]^{vir} \in
A_{0}\left(\overline{M}_{h}(X,\beta_{cs})\right)$$ is in the
$0^{\text{th}}$ Chow group.

Now let $X$ be a $\mathbb{P}^{2}$-bundle of the form
$$\mathbb{P}(L_{0}\oplus L_{1}\oplus L_{2})\rightarrow C,$$ where $C$
is a curve of genus $g$, and $L_{0}\rightarrow C$ is the trivial
bundle, and $L_{1}\rightarrow C$ and $L_{2}\rightarrow C$ are two
arbitrary line bundles of degrees $k_{1}$ and $k_{2}$,
respectively. As in \cite{Bryan-Pandharipande-in-prep}, we use the
word \textbf{level} to refer to the degree of $L_{1}$ or $L_{2}$.
Sometimes we use the notation $\mathcal{O}$ for the trivial bundle
instead of $L_{0}$.

By Leray-Hirsch theorem we have
\[ H^{\text{even}}(X,\mathbb{Z})=
\mathbb{Z}[H,F]/\left(H^{3}+(k_{1}+k_{2})F \cdot H^{2}\right),
\]
where $H$ is the class of the divisor $$\mathbb{P}(L_{1}\oplus
L_{2})\subset \mathbb{P}(L_{0}\oplus L_{1}\oplus L_{2}),$$ and $F$
is the class of the fiber of the bundle $X$. Note also that $H$ is
cohomologous to the first Chern class of the anti-tautological
bundle over $X$.

In this cohomology ring, we have
$$H^{3}=-(k_{1}+k_{2}),$$ $$F \cdot H^{2}=1,$$
and it can be shown that the canonical class of $X$ is given by
\[
K_{X}=-3H+(2g-2-k_{1}-k_{2})F.
\]

\begin{defn}
There is a \textbf{distinguished section} in $X$ which is by
definition the locus of $(1:0:0)$ in
$$X=\mathbb{P}(\mathcal{O}\oplus L_{1}\oplus L_{2})\rightarrow C.$$ We
denote by $\beta_{0}$ the cohomology class in
$H^{4}(X,\mathbb{Z})$ which is represented by this locus. We also
define $$f:=H\cdot F \in H^{4}(X,\mathbb{Z}).$$
\end{defn}
$\{\beta_{0},f\}$ is a set of generators for
$H^{4}(X,\mathbb{Z})$. $\{H^{2},f\}$ is another set of generators.
It is not hard to see that $H \cdot \beta_{0}=0$ and $F \cdot
\beta_{0}=1$, and also
\[
\beta_{0}=H^{2}+ (k_{1}+k_{2})f.
\]

\begin{rem} \label{rem:section class}
One can see easily that for $\mathbb{P}^{2}$-bundles of this form,
$\beta \in H^{4}(X,\mathbb{Z})$ is a section class (see Definition
\ref{defn:section class}) if and only if it is of the form
$$\beta=\beta_{0}+nf$$ for an integer $n$ (see also Remark
\ref{rem:general section class}).
\end{rem}

The complex torus $\mathbb{T}=(\mathbb{C}^{*})^{3}$ acts on
$$X=\mathbb{P}(L_{0}\oplus L_{1}\oplus L_{2})\rightarrow C$$ by
$$(z_{0},z_{1},z_{2})(x_{0}:x_{1}:x_{2}) \mapsto
(z_{0}x_{0}:z_{1}x_{1}:z_{2}x_{2}).$$ Let $\beta_{s} \in
H^{4}(X,\mathbb{Z})$ be a section class. The \textbf{partition
function} of the class $\beta_{s}$ equivariant Gromov-Witten
invariants of $X$ is given by:
\begin{align*}
Z_{\beta_{s}}(g\operatorname{|}k_{1},k_{2})&=\sum_{h=0}^{\infty}u^{2h-2-K_{X}\cdot
\beta_{s}}\int_{[\overline{M}_{h}(X,\beta_{s})]^{vir}}1,
\end{align*}
where $\overline{M}_{h}(X,\beta_{s})$ is the moduli space of
degree $\beta_{s}$, genus $h$ stable maps\footnote{We assume that
all domain curves are connected (see Remark \ref{rem:connected vs
disconnected})} to $X$, and
$$[\overline{M}_{h}(X,\beta_{s})]^{vir}\in A^{\mathbb{T}}_{D}
\left(\overline{M}_{h}(X,\beta_{s})\right)$$ is in the
$D^{\text{th}}$ equivariant Chow group for $$D=-K_{X}\cdot
\beta_{s}=\text{virdim }\overline{M}_{h}(X,\beta_{s}).$$ Since we
are working equivariantly, our definition makes sense even for
negative values of $D$ (c.f. Section 2.2.1 of
\cite{Bryan-Pandharipande-in-prep}).
\begin{rem} \label{rem:equivariant deformation}
The equivariant Gromov-Witten partition functions are invariant
under equivariant deformations. The space
$$X=\mathbb{P}(\mathcal{O}\oplus L_{1}\oplus L_{2})\rightarrow
C$$ that we work with is determined up to equivariant deformation
by $g$, the genus of $C$, and the levels $k_{1}$ and $k_{2}$ of
$L_{1}$ and $L_{2}$, so in this paper we can refer to $X$ by
specifying only these parameters.
\end{rem}

Let $t_{0}, t_{1}, t_{2}$ be the generators for the equivariant
cohomology of a point: \[ H_{\mathbb{T}}^{*}(pt)=
H^{*}((\mathbb{CP}^{\infty})^{3})\cong\mathbb{Q}[t_{0}, t_{1},
t_{2}].
\]
$Z_{\beta_{s}}(g\operatorname{|}k_{1},k_{2})$ is a homogeneous
polynomial in $t_{0}, t_{1}, t_{2}$ of degree $-D$ with
coefficients in $\mathbb{Q}((u))$. In particular, it is zero if
$D$ is positive, and it is a Laurent series in $u$, independent of
$t_{0},t_{1}, t_{2}$, when $D=0$. In the later case,
$Z_{\beta_{s}}(g\operatorname{|}k_{1},k_{2})$ is equal to the
usual Gromov-Witten partition function. (c.f. Section 2.2.1 of
\cite{Bryan-Pandharipande-in-prep}).

The \textbf{partition function} of the section class equivariant
Gromov-Witten invariants of the space $X$ is given by:
\begin{align*}
Z(g\operatorname{|}k_{1},k_{2})=\sum_{\beta_{s} \text{ is a
section class}} Z_{\beta_{s}}(g\operatorname{|}k_{1},k_{2}).
\end{align*}
Note that we can recover any partition function
$Z_{\beta_{s}}(g\operatorname{|}k_{1},k_{2})$ from
$Z(g\operatorname{|}k_{1},k_{2})$ by looking at terms in
$Z(g\operatorname{|}k_{1},k_{2})$ homogeneous in $t_{0},t_{1},
t_{2}$ (see Remark \ref{rem:sum finite}).

In Section \ref{sec:result}, we will prove the main result of this
paper that gives a formula for
$$Z(g\operatorname{|}k_{1},k_{2})$$ for any given genus $g$ and
levels $k_{1}$, $k_{2}$:
\begin{thm} \label{thm:matrix formula}
Let $X$ be a $\mathbb{P}^{2}$ bundle over a curve $C$ of genus $g$
of the form $$\mathbb{P}(\mathcal{O}\oplus\ L_{1} \oplus\
L_{2})\rightarrow C,$$ where $L_{1}$ and $L_{2}$ are two line
bundles of degrees $k_{1}$ and $k_{2}$, respectively, then
$$Z(g\operatorname{|}k_{1},k_{2})=
\operatorname{tr}\left(G^{g-1}U_{1}^{k_{1}}U_{2}^{k_{2}}\right),$$
where the matrices $G$, $U_{1}$ and $U_{2}$ with entries in the
ring $\mathbb{Q}((u))(t_{0},t_{1},t_{2})$ are given by
\begin{align*}
G=&\left[\begin{array}{ccc} (t_{0}-t_{1})(t_{0}-t_{2}) & 0
& 0 \\
0 & (t_{1}-t_{0})(t_{1}-t_{2})
& 0 \\
0 & 0 & (t_{2}-t_{0})(t_{2}-t_{1})
\end{array}
\right] \\
+&\left[\begin{array}{ccc}
\frac{2(2t_{0}-t_{1}-t_{2})}{(t_{0}-t_{1})(t_{0}-t_{2})} &
\frac{t_{0}+t_{1}-2t_{2}}{(t_{0}-t_{1})(t_{0}-t_{2})}
& \frac{t_{0}+t_{2}-2t_{1}}{(t_{0}-t_{1})(t_{0}-t_{2})} \\
\frac{t_{0}+t_{1}-2t_{2}}{(t_{1}-t_{0})(t_{1}-t_{2})} &
\frac{2(2t_{1}-t_{0}-t_{2})}{(t_{1}-t_{0})(t_{1}-t_{2})}
& \frac{t_{1}+t_{2}-2t_{0}}{(t_{1}-t_{0})(t_{1}-t_{2})} \\
\frac{t_{0}+t_{2}-2t_{1}}{(t_{2}-t_{0})(t_{2}-t_{1})} &
\frac{t_{1}+t_{2}-2t_{0}}{(t_{2}-t_{0})(t_{2}-t_{1})} &
\frac{2(2t_{2}+t_{0}-t_{1})}{(t_{2}-t_{0})(t_{2}-t_{1})}
\end{array} \right]\phi^{3},\\
U_{1}=&\left[
\begin{array}{ccc}
\frac{\phi}{t_{0}-t_{1}} &
\frac{(t_{1}-t_{2})\phi}{(t_{0}-t_{1})(t_{0}-t_{2})} & 0 \\
\frac{\phi}{t_{1}-t_{0}} &
+\frac{(t_{1}-t_{0})^{2}(t_{1}-t_{2})^{2}\phi^{-2}+(2t_{1}-t_{0}-t_{2})\phi}
{(t_{1}-t_{0})(t_{1}-t_{2})} & \frac{\phi}{t_{1}-t_{2}} \\
0 & \frac{(t_{1}-t_{0})\phi}{(t_{2}-t_{0})(t_{2}-t_{1})} &
\frac{\phi}{t_{2}-t_{1}}
\end{array} \right],\\
U_{2}=&\left[
\begin{array}{ccc}
\frac{\phi}{t_{0}-t_{2}} & 0 &
\frac{(t_{2}-t_{1})\phi}{(t_{0}-t_{1})(t_{0}-t_{2})} \\
0 & \frac{\phi}{t_{1}-t_{2}} &
\frac{(t_{2}-t_{0})\phi}{(t_{1}-t_{0})(t_{1}-t_{2})} \\
\frac{\phi}{t_{2}-t_{0}} & \frac{\phi}{t_{2}-t_{1}} &
\frac{(t_{2}-t_{0})^{2}(t_{2}-t_{1})^{2}\phi^{-2}+
(2t_{2}-t_{0}-t_{1})\phi}{(t_{2}-t_{0})(t_{2}-t_{1})}
\end{array} \right],
\end{align*}
where $\phi=2\operatorname{sin}\frac{u}{2}$.
\end{thm}
In section \ref{sec:proof of theorem}, as an application of
Theorem \ref{thm:matrix formula}, we will prove the following
result:

\begin{thm} \label{thm:Calabi-Yau class} Let $X\rightarrow C$
be any $\mathbb{P}^{2}$-bundle over a curve $C$ of genus $g$, and
let $\beta_{cs}\in H^{4}(X,\mathbb{Z})$ be a Calabi-Yau section
class, then
\[ Z_{\beta_{cs}}(g)=3^{g}\left({2\operatorname{sin}
\frac{u}{2}}\right)^{2g-2}.
\]
\end{thm}
\subsection*{Plan of the paper}
In Section \ref{sec:relative invariants}, we define the partition
function of the section class relative equivariant Gromov-Witten
invariants of the space
$$X=\mathbb{P}(\mathcal{O}\oplus L_{1}\oplus L_{2})\rightarrow C.$$
Then we express a gluing theorem for these partition functions.

In Section \ref{sec:calculation}, we compute some of the basic
partition functions we defined in Section \ref{sec:relative
invariants}, in the case $g=0$. There are some basic partition
functions in this case that we can compute via localization, we
compute them in \ref{sec:localization}. We use the gluing theorem
of Section \ref{sec:relative invariants} to compute those that we
cannot compute via localizations. This will be done in Section
\ref{subsec:gluing rules}.

In Section \ref{sec:result}, using the results of Section
\ref{sec:calculation}, we construct the matrices $G$, $U_{1}$ and
$U_{2}$ appeared in Theorem \ref{thm:matrix formula} and then we
prove the theorem.

In Section \ref{sec:proof of theorem}, we first prove (Lemma
\ref{lem:deformation}) that any $\mathbb{P}^{2}$-bundle over a
curve $C$ is deformation equivalent to a $\mathbb{P}^{2}$-bundle
over $C$ of the form
$$\mathbb{P}(\mathcal{O} \oplus\mathcal{O} \oplus L).$$
Having this, we use Theorem \ref{thm:matrix formula} to prove
Theorem \ref{thm:Calabi-Yau class}.

In Appendix \ref{app:proof of gluing}, we first prove that it is
enough in this paper to only consider the moduli space of maps
with connected domains (Lemma \ref{lem:connected vs
disconnected}). After that we give a proof for the gluing theorem
expressed in Section \ref{sec:relative invariants}.

In Appendix \ref{app:special cases}, we prove formulas for the
partition function of equivariant Gromov-Witten invariants of the
space $X$ for some special cases of class, $\beta$, and levels,
$k_{1}$ and $k_{2}$.

\subsection*{Acknowledgment} This paper is part of my Ph.D.
thesis at the University of British Columbia. I would like to
express deep gratitude to my supervisor Dr. Jim Bryan whose
guidance was crucial for the successful completion of this work.

\section{Relative invariants and the gluing theorem} \label{sec:relative
invariants} Let $(C,p_{1},\dots,p_{r})$ be a nonsingular curve of
genus $g$ with $r$ marked points. Following the notations of
Section \ref{sec:notation}, we take
$$X=\mathbb{P}(\mathcal{O}\oplus L_{1}\oplus L_{2})\rightarrow
(C,p_{1},\dots,p_{r}).$$ We will review the definition of the
section class equivariant Gromov-Witten invariants
\textbf{relative} to divisors $F_{1},\dots,F_{r}$, where $F_{i}$
is the fiber over the point $p_{i}$. For a treatment of the
foundations of equivariant relative Gromov-Witten theory, see
\cite{Graber-Vakil}.

The complex torus $\mathbb{T}=(\mathbb{C}^{*})^{3}$ acts on $X$ as
in Section \ref{sec:notation}. We need to fix a basis,
$\mathcal{B}_{p}$, for the equivariant cohomology of each fiber,
$F_{p}$, which is a copy of $\mathbb{P}^2$:
$$H^{*}_{\mathbb{T}}(F_{p})\cong
H^{*}_{\mathbb{T}}(\mathbb{P}^{2})\cong
\mathbb{Z}[H](t_{0},t_{1},t_{2})/\left(\prod_{j=0}^{2}(H-t_{j})\right).$$

Let $\beta_{s} \in H^{4}(X,\mathbb{Z})$ be a section class
(defined in Section \ref{sec:notation}). We take
$$Z_{\beta_{s}}^{h}(g\operatorname{|}k_{1},k_{2})_{\alpha_{1}\dots\alpha_{r}}$$ to
be class $\beta_{s}$, genus $h$, equivariant Gromov-Witten
invariant of $X$ relative to the divisors $F_{1},\dots,F_{r}$,
with restrictions given by $\alpha_{p}\in \mathcal{B}_{p}$, one
for each divisor. More precisely, we take
\begin{align*}
&\vec{L}=(l_{1},\dots,l_{r}) \in
\left(\mathbb{Z}^{+}\right)^{r},\\
&\vec{F}=(F_{1},\dots,F_{r}).
\end{align*}
Then following Section 2 of \cite{MNOP2}, let $X[\vec{L}]$ be the
$l_{i}$-step degeneration of $X$ along each $F_{i}$, and let
$$\overline{M}_{h}(X/\vec{F},\beta_{s})$$ be the moduli
space of relative stable maps $$\left[q:C'\rightarrow
X[\vec{L}]\right]$$ from nodal genus $h$ curves\footnote{We assume
that all domain curves are connected (see Remark
\ref{rem:connected vs disconnected})}, $C'$, to $X[\vec{L}]$, for
some $\vec{L}$, which are representing the class $\beta_{s}$. Then
$\overline{M}_{h}(X/\vec{F},\beta_{s})$ is a DM-stack of virtual
dimension $-K_{X}\cdot \beta_{s}$ (see also \cite{Li-relative1}).
\begin{rem}
Since $F \cdot \beta_{s}=1$,
$\overline{M}_{h}(X/\vec{F},\beta_{s})$ does not involve partition
vectors, as it does in a more general case in \cite{MNOP2}. Our
moduli space is more general than the one in \cite{MNOP2} in the
sense that it parameterizes maps relative to more than or equal to
one divisor.
\end{rem}

For each $p$, we have a evaluation map which is determined by
relative points, and is $\mathbb{T}$-equivariant (see
\cite{Li-relative2}):
\[\text{ev}_{p}: \overline{M}_{h}(X/\vec{F},\beta_{s})
\rightarrow F_{p}. \] Then
\[
Z_{\beta_{s}}^{h}(g\operatorname{|}k_{1},k_{2})_{\alpha_{1},\dots,\alpha_{r}}=
\int_{[\overline{M}_{h}(X/\vec{F},\beta_{s})]^{\text{vir}}}
\text{ev}^{*}_{1}(\alpha_{1})\cup \dots \cup
\text{ev}^{*}_{r}(\alpha_{r}),
\]
where $$[\overline{M}_{h}(X,\beta_{s})]^{vir}\in
A^{\mathbb{T}}_{D} \left(\overline{M}_{h}(X,\beta_{s})\right)$$ is
in the $D^{\text{th}}$ equivariant Chow group for
$$D=-K_{X}\cdot \beta_{s}=\text{virdim
}\overline{M}_{h}(X/\vec{F},\beta_{s}).$$ Note that the invariants
can be non-zero even for negative values of $D$ (c.f. Section
2.2.1 of \cite{Bryan-Pandharipande-in-prep}, and also see Remark
\ref{rem:equivariant deformation}). Then the \textbf{partition
function} of the class $\beta_{s}$, relative, equivariant
Gromov-Witten invariants of the space $X$ relative to $\vec{F}$
with relative multiplicities $\alpha_{p}\in \mathcal{B}_{p}$ is
given by:
\begin{equation} \label{equ:class betas inv}
Z_{\beta_{s}}(g\operatorname{|}k_{1},k_{2})_{\alpha_{1}\dots\alpha_{r}}=
\sum_{h=0}^{\infty}
Z_{\beta_{s}}^{h}(g)_{\alpha_{1}\dots\alpha_{r}}
u^{2h-2-K_{X}\cdot\beta_{s}}.
\end{equation}

We can also write the \textbf{partition function} of the section
class, relative, equivariant Gromov-Witten invariants of the space
$X$ relative to $\vec{F}$ with relative multiplicities
$\alpha_{p}\in \mathcal{B}_{p}$ as
\begin{align*}
Z(g\operatorname{|}k_{1},k_{2})_{\alpha_{1}...\alpha_{r}}=\sum_{\beta_{s}
\text{ is a section class}}
Z_{\beta_{s}}(g\operatorname{|}k_{1},k_{2})_{\alpha_{1}...\alpha_{r}}.
\end{align*}
It is evident that when $r=0$ we get the partition function for
the ordinary invariants, defined in Section \ref{sec:notation}.
\begin{rem} \label{rem:global theory}
$Z_{\beta_{s}}(g\operatorname{|}k_{1},k_{2})_{\alpha_{1}...\alpha_{r}}$
is a homogeneous polynomial in $t_{0}, t_{1}, t_{2}$ of degree
$$N=\sum_{p=1}^{r}\operatorname{deg}(\alpha_{p})-D,$$ with
coefficients in $\mathbb{Q}((u))$.\\ In particular it is zero if
$N<0$, and it is a Laurent series in $u$, independent of $t_{0},
t_{1}, t_{2}$, when $N=0$ (c.f.
\cite{Bryan-Pandharipande-in-prep}, Section 2.2.1).
\end{rem}
\begin{rem} \label{rem:sum finite}
We can reexpress the definition of the partition function for the
section class invariants as follows (see Remark \ref{rem:section
class}):
\begin{align*}
Z(g\operatorname{|}k_{1},k_{2})_{\alpha_{1}...\alpha_{r}}=\sum_{n
\in \mathbb{Z}}
Z_{\beta_{0}+nf}(g\operatorname{|}k_{1},k_{2})_{\alpha_{1}...\alpha_{r}}.
\end{align*}
This sum is finite, because by Remark \ref{rem:global theory}, it
is clear that the sum is terminated from above, and it is also
terminated from below because for the large negative values of
$n$, there is no curve representing the class $\beta_{0}+nf$ which
means that
$$\overline{M}_{h}(X/\vec{F},\beta_{0}+nf)=\emptyset$$
for $n\ll0$. To see the last claim, let $E=\mathcal{O}\oplus
L_{1}\oplus L_{2}$, and notice that for a negative $n$ there is a
one to one correspondence between geometric sections representing
$\beta_{0}+nf$ and degree $-n$ sub-line bundles of $E$. $E$ has no
sub-line bundle of degree greater than $k_{1}+k_{2}$. Therefore
for $n\ll0$, $E$ has no sub-line bundle of degree $-n$, which
proves our claim.

One can also recover
$Z_{\beta_{0}+nf}(g\operatorname{|}k_{1},k_{2})_{\alpha_{1}...\alpha_{r}}$
from $Z(g\operatorname{|}k_{1},k_{2})_{\alpha_{1}...\alpha_{r}}$,
by looking at terms in the sum above which are homogeneous in
$t_{0}, t_{1}, t_{2}$ of degree
\begin{align*}
N&=\sum_{p=1}^{r}\operatorname{deg}(\alpha_{p})-D\\
&=\sum_{p=1}^{r}\operatorname{deg}(\alpha_{p})+2g-2-k_{1}-k_{2}-3n.
\end{align*}
\end{rem}

We use the following useful lemma in Section
\ref{sec:calculation}:
\begin{lem} \label{lem:beta=beta+f} Let
$X=\mathbb{P}(\mathcal{O}\oplus L\oplus L')$, where $L$ and $L'$
have levels $n$ and $m$, respectively. Let $\xi_{n}$ be the class
which is represented by the locus of $(0:1:0)$ in
$$X\cong \mathbb{P}(L^{-1}\oplus \mathcal{O} \oplus L'
L^{-1}).$$ Then we have the following relation in
$H^{4}(X,\mathbb{Z})$:
\[
\xi_{n}=\beta_{0}-nf.
\]
\end{lem}
\textsc{Proof.} First notice that there exist $a,b\in\mathbb{Z}$
such that
$$\xi_{n}=a\beta_{0}+bf.$$ By trivial relations
$$F \cdot \beta_{0}=F \cdot \xi_{n}=1$$ and also $F \cdot F=0$,
it is clear that $a=1$. \\For finding $b$, notice that the normal
bundle of $\xi_{n}$ in $X$ is isomorphic to $$L\oplus
L'L^{-1}\rightarrow C,$$ which means that $$H \cdot \xi_{n}=-n$$
(Recall from Section \ref{sec:notation} that $H$ is cohomologous
to $\mathbb{P}(L \oplus L')\subset X$). Now combining this with
$H\cdot \beta_{0}=0$ and $H\cdot f=1$, we see that $b=-n$, which
proves the lemma. \qed

Before expressing the gluing theorem, we fix a basis,
$\mathcal{B}$, for the equivariant cohomology of $\mathbb{P}^{2}$.
We take
\begin{align*}
&x_{0}:=(H-t_{1})(H-t_{2}), \\
&x_{1}:=(H-t_{0})(H-t_{1}), \\
&x_{2}:=(H-t_{0})(H-t_{2}).
\end{align*}
We have $x_{i}\in H^{4}_{\mathbb{T}}(\mathbb{P}^{2})$.
$x_{0},x_{1},x_{2}$ are in fact equivariant cohomology classes
represented by three fixed points of the torus action on
$\mathbb{P}^{2}$. We define
$$\mathcal{B}:=\{x_{0},x_{1},x_{2}\}.$$ It is easy to see that
$\mathcal{B}$ is a basis for
$H^{*}_{\mathbb{T}}(\mathbb{P}^{2})\otimes\mathbb{Q}(t_{0},t_{1},t_{2})$,
for example, we can recover the ordinary basis for the cohomology
of $\mathbb{P}^{2}$:
\begin{align*}
&1\hspace{3mm}=\frac{x_{0}}{(t_{0}-t_{1})(t_{0}-t_{2})}+
\frac{x_{1}}{(t_{1}-t_{0})(t_{1}-t_{2})}+
\frac{x_{2}}{(t_{2}-t_{0})(t_{2}-t_{1})},\\
&H\hspace{1.5mm}=\frac{t_{0}x_{0}}{(t_{0}-t_{1})(t_{0}-t_{2})}+
\frac{t_{1}x_{1}}{(t_{1}-t_{0})(t_{1}-t_{2})}+
\frac{t_{2}x_{2}}{(t_{2}-t_{0})(t_{2}-t_{1})},\\
&H^{2}=\frac{t_{0}^{2}x_{0}}{(t_{0}-t_{1})(t_{0}-t_{2})}+
\frac{t_{1}^{2}x_{1}}{(t_{1}-t_{0})(t_{1}-t_{2})}+
\frac{t_{2}^{2}x_{2}}{(t_{2}-t_{0})(t_{2}-t_{1})}.
\end{align*}

We also have these relations:
\begin{align} \label{equ:basis relations}
x_{0}^{2}&=(t_{0}-t_{1})(t_{0}-t_{2})x_{0},\nonumber\\
x_{1}^{2}&=(t_{1}-t_{0})(t_{1}-t_{2})x_{1},\nonumber\\
x_{2}^{2}&=(t_{2}-t_{0})(t_{2}-t_{1})x_{2},
\nonumber\\
x_{i}x_{j}&=0 \text{\hspace{.5cm}for\hspace{.15cm}} i\neq j.
\end{align}

\textbf{Convention.} From now on, we assume that each $\alpha_{p}$
for $p=1,\dots,r$, in the definition of relative partition
functions belong to this basis set, $\mathcal{B}$, with the
identification of each $F_{p}$ with $\mathbb{P}^{2}$.

We take
\begin{align*}
&T(x_{0}):=(t_{0}-t_{1})(t_{0}-t_{2}),\\
&T(x_{1}):=(t_{1}-t_{0})(t_{1}-t_{2}),\\
&T(x_{2}):=(t_{2}-t_{0})(t_{2}-t_{1}).
\end{align*}
Then we raise the indices for the relative partition functions by
the following rule:
$$Z(g\operatorname{|}k_{1},k_{2})_{\alpha_{1}\dots\alpha_{s}}
^{\gamma_{1}\dots\gamma_{t}}:=\left(\prod_{p=1}^{t}\frac{1}{T(\gamma_{p})}\right)
Z(g\operatorname{|}k_{1},k_{2})_{\alpha_{1}\dots\alpha_{s}
\gamma_{1}\dots\gamma_{t}}.$$

Then we have the following gluing rules similar to Theorem 3.1 in
\cite{Bryan-Pandharipande-in-prep}:
\begin{thm} \label{thm:gluing rules}
For any choices of elements $\alpha_{1},\dots,\alpha_{s}$ and
$\gamma_{1},\dots,\gamma_{t}$ from the set $\mathcal{B}$, and
integers satisfying $g=g'+g''$, $k_{1}=k_{1}'+k_{1}''$, and
$k_{2}=k_{2}'+k_{2}''$ we have
\begin{align*}
Z(g\operatorname{|}k_{1},k_{2})_{\alpha_{1}\dots\alpha_{s}
\gamma_{1}\dots\gamma_{t}}=\sum_{\lambda \in \mathcal{B}}
Z(g'\operatorname{|}k_{1}',k_{2}')_{\alpha_{1}...\alpha_{s}\lambda}
Z(g''\operatorname{|}k_{2}'',k_{2}'')_{\gamma_{1}\dots\gamma_{t}}^{\lambda},
\end{align*}
and
\[Z(g\operatorname{|}k_{1},k_{2})_{\alpha_{1}\dots\alpha_{s}}=
\sum_{\lambda \in \mathcal{B}}
Z(g\operatorname{|}k_{1},k_{2})_{\alpha_{1}\dots\alpha_{s}\lambda}^{\lambda}.\]
\end{thm}
The proof of this theorem will be given in Appendix \ref{app:proof
of gluing}.
\begin{rem} \label{rem:connected vs disconnected}
In most of the contexts in which the relative Gromov-Witten
invariants are being used, maps with disconnected domain curves
are considered as well as ones with connected domains. In Lemma
\ref{lem:connected vs disconnected} we prove that in our case,
where we only deal with section classes, we don't need to consider
disconnected domain curves.
\end{rem}
\begin{rem} \label{rem:TQFT}
In exactly the same way as in \cite{Bryan-Pandharipande-in-prep},
one can prove by using Theorem \ref{thm:gluing rules} that the
partition functions
$Z(g\operatorname{|}0,0)_{\alpha_{1}\dots\alpha_{r}}$ give rise to
a $1+1$-dimensional TQFT taking values in the ring
$R=\mathbb{Q}((u))(t_{0},t_{1},t_{2})$. The Frobenius algebra
corresponding to this TQFT (see \cite{Br-Pa-TQFT}, Theorem 2.1) is
$$H=\bigoplus_{i=0}^{2}Re_{x_{i}}$$ for $x_{i}\in \mathcal{B}$,
with multiplication given by
$$
e_{x_{i}}\otimes
e_{x_{j}}=\sum_{k=0}^{2}Z(g\operatorname{|}0,0)_{x_{i}x_{j}}^{x_{k}}e_{k}.$$
We will prove that this Frobenius algebra and hence the
corresponding TQFT is semisimple (Corollary \ref{pro:semisimple}).
In Section \ref{sec:proof of theorem}, we use this fact for
proving Theorem \ref{thm:matrix formula} and \ref{thm:Calabi-Yau
class} and also the results in Appendix \ref{app:special cases},
for the case $g=0$.
\end{rem}

We will use the following corollary of Theorem \ref{thm:gluing
rules} in our calculations:
\begin{cor} \label{cor:gluing rule}
With the same notation as in Theorem \ref{thm:gluing rules} we
have
\begin{align*}
Z_{\beta_{0}+nf}&(g\operatorname{|}k_{1},k_{2})_{\alpha_{1}...\alpha_{s}
\gamma_{1}\dots\gamma_{t}}=\\&\sum_{\lambda \in \mathcal{B}}
\sum_{n=n'+n''}
Z_{\beta_{0}+n'f}(g'\operatorname{|}k_{1}',k_{2}')_{\alpha_{1}...\alpha_{s}\lambda}
Z_{\beta_{0}+n''f}(g''\operatorname{|}k_{2}'',k_{2}'')
_{\gamma_{1}\dots\gamma_{t}}^{\lambda}.
\end{align*}
\end{cor}\qed

\section{Calculations} \label{sec:calculation}
We work with the space
$$X=\mathbb{P}(\mathcal{O}\oplus L_{1}\oplus L_{2})\rightarrow
(C,p_{1},\dots,p_{r}),$$ throughout this section. In accordance
with the notations in \cite{Bryan-Pandharipande-in-prep}, we will
use the words \textbf{cap}, \textbf{tube} and \textbf{pants} to
refer to the case where the base curve, $C$, is a genus zero curve
with one, two and three marked points, respectively, and by
$(k_{1},k_{2})$ we mean that the level of the line bundle $L_{i}$
is $k_{i}$ for $i=1,2$ (see Remark \ref{rem:equivariant
deformation}). We sometimes refer to the partition functions by
referring to the space to which they correspond. Finally, for
simplicity, we will use the notation
$$\phi:=2\operatorname{sin}\frac{u}{2}$$ in later calculations.

Similar to Section $4.3$ in \cite{Bryan-Pandharipande-in-prep} one
can see that the following partition functions determine the
theory completely:
\begin{eqnarray*}
Z(0\operatorname{|}0,0)_{\alpha}\quad &:& \quad
\text{corresponding to the level $(0,0)$ cap},
\\Z(0\operatorname{|}0,0)_{\alpha_{1}\alpha_{2}} \quad &:& \quad
\text{corresponding to the level $(0,0)$ tube},
\\Z(0\operatorname{|}0,0)_{\alpha_{1}\alpha_{2}\alpha_{3}}\quad &:& \quad
\text{corresponding to the level $(0,0)$ pants},
\\Z(0\operatorname{|}-1,0)_{\alpha}\quad &:& \quad
\text{corresponding to the level $(-1,0)$ cap},
\\Z(0\operatorname{|}0,-1)_{\alpha}\quad &:& \quad \text{corresponding to the level
$(0,-1)$ cap},
\\Z(0\operatorname{|}1,0)_{\alpha}\quad &:& \quad\text{corresponding to the level
$(1,0)$ cap},
\\Z(0\operatorname{|}0,1)_{\alpha}\quad &:& \quad\text{correspond to the level
$(0,1)$ cap}.
\end{eqnarray*}
We refer to the partition functions above as the \textbf{basic}
partition functions.

By the discussion given in Remark \ref{rem:sum finite}, one can
prove the following lemma:
\begin{lem} \label{lem:tqft pieces}
The basic partition functions are given by
\begin{eqnarray*}
Z(0\operatorname{|}0,0)_{\alpha}\quad&=&\quad
Z_{\beta_{0}}(0\operatorname{|}0,0)_{\alpha},
\\Z(0\operatorname{|}0,0)_{\alpha_{1}\alpha_{2}}\quad&=&\quad
Z_{\beta_{0}}(0\operatorname{|}0,0)_{\alpha_{1}\alpha_{2}},
\\Z(0\operatorname{|}0,0)_{\alpha_{1}\alpha_{2}\alpha_{3}}\quad&=&\quad
Z_{\beta_{0}}(0\operatorname{|}0,0)_{\alpha_{1}\alpha_{2}\alpha_{3}}+
Z_{\beta_{0}+f}(0\operatorname{|}0,0)_{\alpha_{1}\alpha_{2}\alpha_{3}},
\\Z(0\operatorname{|}-1,0)_{\alpha}\quad&=&\quad
Z_{\beta_{0}}(0\operatorname{|}-1,0)_{\alpha},
\\Z(0\operatorname{|}0,-1)_{\alpha} \quad&=&\quad
Z_{\beta_{0}}(0\operatorname{|}0,-1)_{\alpha},
\\Z(0\operatorname{|}1,0)_{\alpha}\quad&=&\quad
Z_{\beta_{0}-f}(0\operatorname{|}1,0)_{\alpha},
\\Z(0\operatorname{|}0,1)_{\alpha} \quad&=
&\quad Z_{\beta_{0}-f}(0\operatorname{|}0,1)_{\alpha}.
\end{eqnarray*}
\end{lem}
\textsc{Proof:}
We prove the third equality as follows: \\
In the right hand side, we do not have any partition function of
class $\beta_{0}+nf$ for $n<0$, because $\mathcal{O}\oplus
L_{1}\oplus L_{2}$ does not have any sub-line bundle of a positive
degree, as $L_{1}$ and $L_{2}$ are level zero.

We also do not have any partition function of class $\beta_{0}+nf$
for $n>1$, because
\begin{align*}
N&=\sum_{p=1}^{3}\operatorname{deg}(\alpha_{p})-D\\
&=(2+2+2)-(3H+2F)\cdot(\beta_{0}+nf)\\&=4-3n
\end{align*}
which is negative for $n>1$ (see Remark \ref{rem:global
theory}).\\ The other equalities are proved similarly. \qed

The rest of this section is devoted to computing the terms
appeared in the right hand sides of the equations in Lemma
\ref{lem:tqft pieces}.

\subsection{Calculations via localization} \label{sec:localization}
The complex torus acts on $X$ as before. We define
\begin{align*}
&S_{0}:\quad \text{the locus of \hspace{1mm}$(1:0:0)$\hspace{1mm}
in \hspace{1mm}} X\cong\mathbb{P}(\mathcal{O}\oplus
L_{1}L_{0}^{-1}\oplus L_{2}L_{0}^{-1}),
\\&S_{1}:\quad \text{the locus of \hspace{1mm}$(0:1:0)$\hspace{1mm} in
\hspace{1mm}} X\cong \mathbb{P}(L_{0}L_{1}^{-1}\oplus
\mathcal{O}\oplus L_{2} L_{1}^{-1}),
\\&S_{2}: \quad \text{the locus of \hspace{1mm}$(0:0:1)$\hspace{1mm} in
\hspace{1mm}} X\cong \mathbb{P}(L_{0}L_{2}^{-1}\oplus
L_{1}L_{2}^{-1}\oplus \mathcal{O}).
\end{align*}
$S_{0}$, $S_{1}$ and $S_{2}$ are fixed under the torus action, and
by Lemma \ref{lem:beta=beta+f}, they represent the classes
$\beta_{0}$, $\beta_{0}-k_{1}f$ and $\beta_{0}-k_{2}f$,
respectively.

As before, let $\beta_{s}$ be a section class. The torus action on
$X$ induces an action on $\overline{M}_{h}(X/\vec{F},\beta_{s})$.
We denote the fixed locus of this action by
$\overline{M}_{h}(X/\vec{F},\beta_{s})^{\mathbb{T}}.$

By notations of Section \ref{sec:relative invariants}, we let
$S_{i}[\vec{L}]\subset X[\vec{L}]$ be the $l_{i}$-step
degeneration of $S_{i}$ along the intersection point
$$\tau_{ip}=S_{i}\cap F_{p}$$ for $p=1,\dots,r$, and $i=0,1,2$.\\
Then $\overline{M}_{h}(X/\vec{F},\beta_{s})^{\mathbb{T}}$
parameterizes maps $$\left[q:C'\rightarrow X[\vec{L}]\right]$$ for
some $\vec{L}$, whose images are either of
$$S_{i}[\vec{L}]\cup_{n=1}^{m_{i}}f_{n}$$ for $i=0,1$ or $2$,
where by the last expression we mean $S_{i}[\vec{L}]$ with $m_{i}$
$\mathbb{T}$-fixed fiber curves, $f_{n}$ ($f_{n}$ represents the
class $bf$ for some $a \in \mathbb{Z}^{+}$), are attached to it at
some points. Note that the choice of $i \in \{0,1,2\}$, and also
the number of fibers which are attached to $S_{i}[\vec{L}]$,
$m_{i}$, are constrained by the class $\beta_{s}$.

In general, the moduli space
$\overline{M}_{h}(X/\vec{F},\beta_{s})^{\mathbb{T}}$ can be quiet
complicated, because of the existence of the fibers attached to
each $S_{i}[\vec{L}]$. However, in the special case where
$m_{i}=0$ for some $i$, an elementary observation shows that the
component of $\overline{M}_{h}(X/\vec{F},\beta_{s})^{\mathbb{T}}$,
parameterizing maps with images equal to $S_{i}[\vec{L}]$ is
exactly the moduli space of degree one relative stable maps to
curves, which we denote by
$\overline{M}_{h}(S_{i}/\vec{\tau}_{i},1)$, where
$$\vec{\tau}_{i}=(\tau_{i1},\dots,\tau_{ip}).$$
\textbf{Assumption 1.} For the rest of Section
\ref{sec:localization}, we assume that
\begin{equation} \label{equ:easy case for locln}
\overline{M}_{h}(X/\vec{F},\beta_{s})^{\mathbb{T}}= \bigcup_{i\in
I}\overline{M}_{h}(S_{i}/\vec{\tau}_{i},1), \end{equation} where
$I\subset \{0,1,2\}$, depending on the class $\beta_{s}$.

Then one can see that the $\mathbb{T}$-fixed part of the perfect
obstruction theory of $\overline{M}_{h}(X/\vec{F},\beta_{s})$ is
exactly the usual obstruction theory of $\bigcup_{i\in
I}\overline{M}_{h}(S_{i}/\vec{\tau}_{i},1)$, and therefore
$$[\overline{M}_{h}(X/\vec{F},\beta_{s})^{\mathbb{T}}]^{vir}\cong
\sum_{i\in I}[\overline{M}_{h}(S_{i}/\vec{\tau}_{i},1)]^{vir}.$$

In the special case where $m_{i}=0$ for all possible $i$,
Assumption 1 holds. One can see easily that this is the case for
all the partition functions in the right hand sides of equations
in Lemma \ref{lem:tqft pieces}, except for
$$Z_{\beta_{0}+f}(0\operatorname{|}0,0)_{\alpha_{1}\alpha_{2}\alpha_{3}}.$$
These partition functions are calculated in this section via
localization.\\
$Z_{\beta_{0}+f}(0\operatorname{|}0,0)_{\alpha_{1}\alpha_{2}\alpha_{3}}$
will be calculated in Section \ref{subsec:gluing rules} by
combining the results of this section with the results of the
gluing techniques.

Applying the relative virtual localization formula (see Section 3
of \cite{Graber-Vakil}\footnote{In \cite{Graber-Vakil}, the
authors assume for convenience that the relative divisor is in the
fixed locus, but it is straight forward to adapt their methods to
the case at hand \cite{Graber-communication}.} to (\ref{equ:class
betas inv}), we can write
\begin{align} \label{equ:loc res inv}
Z_{\beta_{s}}(&g\operatorname{|}k_{1},k_{2})_{\alpha_{1}\dots
\alpha_{r}}=\nonumber
\\&\sum_{h=0}^{\infty}u^{2h-2-K_{X}\cdot\beta_{s}}
\int_{[\overline{M}_{h}(X/\vec{F},\beta_{s})^{\mathbb{T}}]^{\text{vir}}}
\frac{\text{ev}_{1}^{*}(\alpha_{1})\cap \dots\cap
\text{ev}_{r}^{*}(\alpha_{r})}{e(\text{Norm}^{vir})},
\end{align} where $\text{Norm}^{vir}$ is the equivariant virtual
normal bundle of
$$\overline{M}_{h}(X/\vec{F},\beta_{s})^{\mathbb{T}}\subset
\overline{M}_{h}(X/\vec{F},\beta_{s}),$$ and
$e(\text{Norm}^{vir})$ is its equivariant Euler class.

Let
$$\pi:U\rightarrow\overline{M}_{h}(X/\vec{F},\beta_{s})$$
and
$$q:U\rightarrow \mathcal{X}$$ be the universal curve and universal
map, respectively, where
$$\mathcal{X}\rightarrow\overline{M}_{h}(X/\vec{F},\beta_{s})$$
is the universal target space.

Now notice that the normal bundle of each $F_{p}$ in $X$ is the
trivial bundle with the trivial torus action. So the deformations
of the singularities of the degenerated target spaces does not
contribute in $e(\text{Norm}^{vir})$ (see Section 3 of
\cite{Graber-Vakil}). We also have the following short exact
sequence
$$0\rightarrow \Omega_{X}\rightarrow\Omega_{X}(\operatorname{log}
F_{p})\rightarrow N_{F_{p}|X}\rightarrow 0,$$ where
$\Omega_{X}(\operatorname{log}F_{p})$ is the sheaf of K\"{a}hler
differential with logarithmic poles along $F_{p}$ (see
\cite{Graber-Vakil}), and $N_{F_{p}|X}$ is the normal bundle of
$F_{p}$ in $X$. Again since the torus action on $N_{F_{p}|X}$ is
trivial, from this sequence we can see that the moving part of
$H^{\bullet}(C',q^{*}T_{X}(-\operatorname{log} F_{p}))$ is equal
to the moving part of $H^{\bullet}(C',q^{*}T_{X})$, where
$T_{X}(-\operatorname{log}F_{p})$ is the dual sheaf of
$\Omega_{X}(\operatorname{log}F_{p})$. So from (2) and (3) in
\cite{Graber-Vakil} (the third term in (3) has no moving part), we
can write
\begin{align*}
e(\text{Norm}^{vir})&=e\left(R^{\bullet}\pi_{*}q^{*}(T_{X})^{mov}\right)\\&=\sum_{i\in
I}e\left(R^{\bullet}\pi_{*}q^{*}(N_{S_{i}|X})\right),
\end{align*}
where $N_{S_{i}|X}$ is the normal bundle of $S_{i}$ in $X$. One
can see easily that
\begin{align} \label{equ:normal bundles}
&N_{S_{0}|X}\cong L_{1}L_{0}^{-1}\oplus
L_{2}L_{0}^{-1},\nonumber \\
&N_{S_{1}|X}\cong L_{0}L_{1}^{-1}\oplus
L_{2}L_{1}^{-1},\nonumber \\
&N_{S_{2}|X}\cong L_{0}L_{2}^{-1}\oplus L_{1}L_{2}^{-1}.
\end{align}

For manipulating the evaluation functions in (\ref{equ:loc res
inv}), we use the following cartesian diagram for each
$p=1,\dots,r$, and $i\in I$:
\[
\begin{CD}
\overline{M}_{h}(S_{i}/\vec{\tau}_{i},1) @>>> \{\tau_{ip}\}
\\
@Vj_{i}VV                                       @VVV \\
\overline{M}_{h}(X/\vec{F},\beta_{s})
@>{\text{ev}_{p}}>>F_{p} \\
\end{CD}
\]
where two vertical maps are inclusions, and $\tau_{ip}$ is the
intersection point of $S_{i}$ with $F_{p}$, which is the fixed
point of the torus action on $F_{p}$ representing the class
$x_{i}\in \mathcal{B}$. From this diagram it is clear that
$\text{ev}_{p}^{*}(\alpha_{p})$, restricted to
$\overline{M}_{h}(S_{i}/\vec{\tau}_{i},1)$, is a class of pure
weight for each $p$, and can be taken out of the integrals.

We summarize all the discussion above in the following equation:
\begin{align}\label{equ:final loc res inv}
&Z_{\beta}(g\operatorname{|}k_{1},k_{2})_{\alpha_{1}\dots
\alpha_{r}}=\nonumber
\\&\sum_{h=0}^{\infty}u^{2h-2-K_{X}\cdot\beta_{s}} \sum_{i\in
I}\left(\prod_{p=1}^{r}(\operatorname{ev}_{p}\circ
j_{i})^{*}(\alpha_{p})\right)
\int_{[\overline{M}_{h}(S_{i}/\vec{\tau}_{i},1)]^{\text{vir}}}
e\left(-R^{\bullet}\pi_{*}q^{*}(N_{S_{i}|X})\right).
\end{align}
We will handle the integrals of this form in our calculations by
appealing to the results in \cite{Bryan-Pandharipande-in-prep}.

Applying Atiah-Bott localization theorem to $\mathbb{P}^2$, one
can see easily that
\begin{equation} \label{equ:evaluation}
(\text{ev}_{p}\circ j_{i})^{*}(x_{k})=\begin{cases} T(x_{i}) &
\text{if} \hspace{10pt} i=k \\ 0 & \text{otherwise}
\end{cases} \end{equation} for $k \in \{0,1,2\}$ (see Section
\ref{sec:relative invariants} for the definition of $T(-)$).

\subsubsection{Computing class
$\beta_{0}$, level $(0,0)$ cap, tube and pants:}
\begin{lemma} \label{lemma:level 0 cap}
Partition functions for the class $\beta_{0}$, level $(0,0)$ cap,
tube and pants are given by
\begin{align*}
&Z_{\beta_{0}}(0\operatorname{|}0,0)_{x_{a}}\text{\hspace{6mm}}=1
\\&Z_{\beta_{0}}(0\operatorname{|}0,0)_{x_{a}x_{b}}\text{\hspace{3mm}}=\begin{cases}
T(x_{a}) & \text{if} \hspace{10pt} a=b \\ 0 & \text{otherwise}
\end{cases}
\\&Z_{\beta_{0}}(0\operatorname{|}0,0)_{x_{a}x_{b}x_{c}}=
\begin{cases} T(x_{a})^{2} & \text{if}
\hspace{10pt} a=b=c \\ 0 & \text{otherwise}
\end{cases}
\end{align*}
for $a,b,c \in \{0,1,2\}$.
\end{lemma}
\textsc{Proof:} Since $k_{1}=k_{2}=0$, by Lemma
\ref{lem:beta=beta+f}, all $S_{0}$, $S_{1}$ and $S_{2}$ represent
the class $\beta_{0}$, so in (\ref{equ:final loc res inv}) we have
$I=\{0,1,2\}$.

We use the results of Sections 6.2 and 6.4.2 in
\cite{Bryan-Pandharipande-in-prep} to evaluate the integrals in
(\ref{equ:final loc res inv}), for the cap, tube and pants. We
prove the formula for the tube, the other cases are similar. By
(\ref{equ:final loc res inv}) and Lemma 6.1 in
\cite{Bryan-Pandharipande-in-prep} (for $d=1$) and also
(\ref{equ:normal bundles}) and (\ref{equ:evaluation}) we have
\begin{align*}
&Z_{\beta_{0}}(0\operatorname{|}0,0)_{x_{a}x_{b}}
=\sum_{h=0}^{\infty}u^{2h-2-K_{X}\cdot\beta}\\
\bigg(&\hspace{4.5mm}(\text{ev}_{1}\circ
j_{0})^{*}(x_{a})(\text{ev}_{2}\circ
j_{0})^{*}(x_{b})\\
&\hspace{3.5cm}\cdot\int_{\left[\overline{M}_{h}
(S_{0}/(\tau_{01},\tau_{02}),1)\right]^{\text{vir}}}
e\left(-R^{\bullet}\pi_{*}q^{*}(L_{1}L_{0}^{-1}\oplus
L_{2}L_{0}^{-1})\right)\\
&+(\text{ev}_{1}\circ j_{1})^{*}(x_{a})(\text{ev}_{2}\circ
j_{1})^{*}(x_{b})\\
&\hspace{3.5cm}\cdot
\int_{\left[\overline{M}_{h}(S_{1}/(\tau_{11},\tau_{12}),1)\right]^{\text{vir}}}
e\left(-R^{\bullet}\pi_{*}q^{*}(L_{0}L_{1}^{-1}\oplus
L_{2}L_{1}^{-1})\right)\\
&+(\text{ev}_{1}\circ j_{2})^{*}(x_{a})(\text{ev}_{2}\circ
j_{2})^{*}(x_{b})\\ &\hspace{3.5cm}\cdot
\int_{\left[\overline{M}_{h}(S_{2}/(\tau_{21},\tau_{22}),1)\right]^{\text{vir}}}
e\left(-R^{\bullet}\pi_{*}q^{*}(L_{0}L_{2}^{-1}\oplus
L_{1}L_{2}^{-1})\right)\bigg)\\
&\hspace{2.5cm}=(\delta^{0}_{a}T(x_{0}))(\delta^{0}_{b}T(x_{0}))\frac{1}{T(x_{0})}\\
&\hspace{2.55cm}+
(\delta^{1}_{a}T(x_{1}))(\delta^{1}_{b}T(x_{1}))\frac{1}{T(x_{1})}\\
&\hspace{2.55cm}+(\delta^{2}_{a}T(x_{2}))(\delta^{2}_{b}T(x_{2}))\frac{1}{T(x_{2})}\\
&\hspace{2.5cm}=\delta^{b}_{a}T(x_{a}).\text{\hspace{7.5cm}}\qed
\end{align*}
\subsubsection{Computing class $\beta_{0}$, level
$(0,-1)$ and $(-1,0)$ and class $\beta_{0}-f$, level $(0,1)$ and
$(1,0)$ caps:}
\begin{lemma} \label{lemma:level -1 cap}
Partition functions for the class $\beta_{0}$, level $(0,-1)$ and
$(-1,0)$ caps are given by
\begin{align*}
&Z_{\beta_{0}}(0\operatorname{|}0,-1)_{x_{a}}=
(t_{a}-t_{2})\phi^{-1}
\\&Z_{\beta_{0}}(0\operatorname{|}-1,0)_{x_{a}}=
(t_{a}-t_{1})\phi^{-1}
\end{align*}
for $a=0,1,2$.
\end{lemma}
\textsc{Proof:} We prove the first formula, the second one is
proved in a similar way. We have $k_{1}=0$ and $k_{2}=-1$, so by
Lemma \ref{lem:beta=beta+f}, $S_{0}$, $S_{1}$ represent the class
$\beta_{0}$, but $S_{2}$ represents the class $\beta_{0}+f$.
Therefore in (\ref{equ:final loc res inv}) we have $I=\{0,1\}$.

By Lemma 6.3 in \cite{Bryan-Pandharipande-in-prep} (for $d=1$) and
also (\ref{equ:normal bundles}) and (\ref{equ:evaluation}) we can
rewrite (\ref{equ:final loc res inv}) as
\begin{align*}
Z_{\beta_{0}}(0\operatorname{|}0,-1)_{x_{a}}
&=\sum_{h=0}^{\infty}u^{2h-2-K_{X}\cdot\beta_{0}}\\
&\hspace{-1cm}\bigg(\hspace{2.4mm}(\text{ev}\circ
j_{0})^{*}(x_{a})\int_{[\overline{M}_{h}(S_{0}/\tau_{01},1)]^{\text{vir}}}
e\left(-R^{\bullet}\pi_{*}q^{*}(L_{1}L_{0}^{-1}\oplus
L_{2}L_{0}^{-1})\right)\\&\hspace{-.9cm}+(\text{ev}\circ
j_{1})^{*}(x_{a})\int_{[\overline{M}_{h}(S_{1}/\tau_{11},1)]^{\text{vir}}}
e\left(-R^{\bullet}\pi_{*}q^{*}(L_{0}L_{1}^{-1}\oplus
L_{2}L_{1}^{-1})\right)\bigg)\\
&=\left(\frac{\delta^{0}_{a}T(x_{0})}{t_{0}-t_{1}}+
\frac{\delta^{1}_{a}T(x_{1})}{t_{1}-t_{0}}\right)\phi^{-1}\\
&=(t_{a}-t_{2})\phi^{-1}. \text{\hspace{7.5cm}}\qed
\end{align*}

\begin{lemma} \label{lemma:level 1 cap}
Partition functions for the class $\beta_{0}-f$, level $(0,1)$ and
$(1,0)$ caps are given by
\begin{align*}
&Z_{\beta_{0}-f}(0\operatorname{|}0,1)_{x_{a}}=
(t_{a}-t_{0})(t_{a}-t_{1})\phi^{-2}
\\&Z_{\beta_{0}-f}(0\operatorname{|}1,0)_{x_{a}}=
(t_{a}-t_{0})(t_{a}-t_{2})\phi^{-2}
\end{align*}
for $a=0,1,2$.
\end{lemma}
\textsc{Proof:} We again prove the first relation. The proof of
the second one is similar. In this case only $S_{2}$ represents
the class $\beta_{0}-f$, so we have $I=\{2\}$. The integral that
we need to know in order to prove the Lemma can be found in
Section 8 of \cite{Bryan-Pandharipande-in-prep}. The rest of the
proof is quiet similar to the proof of Lemma \ref{lemma:level -1
cap}.\qed

\subsection{Calculations via gluing techniques}\label{subsec:gluing
rules} In this section, we use Corollary \ref{cor:gluing rule}
(which is referred to as \textbf{the gluing formula}), and the
results of Section \ref{sec:localization} to find
$Z_{\beta_{0}+f}(0\operatorname{|}0,0)_{\alpha_{1}\alpha_{2}\alpha_{3}}$.
For a treatment of gluing spaces and applying the gluing theorem
see Appendix \ref{app:proof of gluing}.

We first need to find the following partition functions of tubes:
\begin{align}\label{equ:tubes}
&Z(0\operatorname{|}0,-1)_{\alpha_{1}\alpha_{2}}=
Z_{\beta_{0}}(0\operatorname{|}0,-1)_{\alpha_{1}\alpha_{2}}\text{\hspace{.8mm}}+
Z_{\beta_{0}+f}(0\operatorname{|}0,-1)_{\alpha_{1}\alpha_{2}},\nonumber
\\&Z(0\operatorname{|}-1,0)_{\alpha_{1}\alpha_{2}}=
Z_{\beta_{0}}(0\operatorname{|}-1,0)_{\alpha_{1}\alpha_{2}}\text{\hspace{.8mm}}+
Z_{\beta_{0}+f}(0\operatorname{|}-1,0)_{\alpha_{1}\alpha_{2}},\nonumber
\\&Z(0\operatorname{|}0,1)_{\alpha_{1}\alpha_{2}}\text{\hspace{3.2mm}}=
Z_{\beta_{0}-f}(0\operatorname{|}0,1)_{\alpha_{1}\alpha_{2}}+
Z_{\beta_{0}}(0\operatorname{|}0,1)_{\alpha_{1}\alpha_{2}},\nonumber
\\&Z(0\operatorname{|}1,0)_{\alpha_{1}\alpha_{2}}\text{\hspace{3.2mm}}=
Z_{\beta_{0}-f}(0\operatorname{|}1,0)_{\alpha_{1}\alpha_{2}}+
Z_{\beta_{0}}(0\operatorname{|}1,0)_{\alpha_{1}\alpha_{2}}.
\end{align}
These equalities can be proved similar to the proof of Lemma
\ref{lem:tqft pieces}. Now we are going to find the partition
functions in the right hand sides of equations in
(\ref{equ:tubes}):

\subsubsection{Computing class $\beta_{0}$, level
$(0,-1)$ and $(-1,0)$ and class $\beta_{0}-f$, level $(0,1)$ and
$(1,0)$ tubes:}
\begin{lemma} \label{lemma:level -1 tube}
Partition functions for the class $\beta_{0}$, level $(0,-1)$ and
$(-1,0)$ tubes are given by
\begin{align*}
&Z_{\beta_{0}}(0\operatorname{|}0,-1)_{x_{a}x_{b}}=
\begin{cases} (t_{0}-t_{1})(t_{0}-t_{2})^{2}\phi^{-1} & \text{if}
\hspace{10pt}
a=b=0,\\
(t_{1}-t_{0})(t_{1}-t_{2})^{2}\phi^{-1} & \text{if} \hspace{10pt}
a=b=1,\\
0 & \text{otherwise.}
\end{cases}
\end{align*}
\begin{align*} &Z_{\beta_{0}}(0\operatorname{|}-1,0)_{x_{a}x_{b}}=
\begin{cases} (t_{0}-t_{2})(t_{0}-t_{1})^{2}\phi^{-1} & \text{if}
\hspace{10pt}
a=b=0,\\
(t_{2}-t_{0})(t_{2}-t_{1})^{2}\phi^{-1} & \text{if} \hspace{10pt}
a=b=2,\\
0 & \text{otherwise}
\end{cases}
\end{align*}for $a,b \in \{0,1,2\}.$
\end{lemma}
\textsc{Proof:} The first relation is simply proved by attaching
the level $(0,-1)$ cap to the level $(0,0)$ pants and applying the
gluing formula. This is schematically indicated by the following
picture:
\begin{center}
\includegraphics{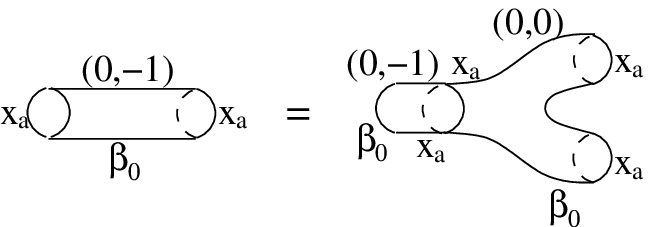}
\end{center}
$$Z_{\beta_{0}}(0\operatorname{|}0,-1)_{x_{a}x_{a}}=
Z_{\beta_{0}}(0\operatorname{|}0,-1)_{x_{a}}
Z_{\beta_{0}}(0\operatorname{|}0,0)_{x_{a}x_{a}}^{x_{a}}.$$ The
result is now obvious by applying Lemma \ref{lemma:level -1 cap}
and \ref{lemma:level 0 cap}. The proof of the second relation is
similar. \qed

Similar to the proof of Lemma \ref{lemma:level -1 tube}, we can
prove this:
\begin{lemma} \label{lemma:level 1 tube}
Partition functions for the class $\beta_{0}-f$, level $(0,1)$ and
$(1,0)$ tubes are given by
\begin{align*}
&Z_{\beta_{0}-f}(0\operatorname{|}0,1)_{x_{a}x_{b}}=
\begin{cases} (t_{2}-t_{0})^{2}(t_{2}-t_{1})^{2}\phi^{-2} & \text{if}
\hspace{10pt}
a=b=2,\\
0 & \text{otherwise.}
\end{cases}
\end{align*}
\begin{align*} &Z_{\beta_{0}-f}(0\operatorname{|}1,0)_{x_{a}x_{b}}=
\begin{cases} (t_{1}-t_{2})^{2}(t_{1}-t_{0})^{2}\phi^{-2} & \text{if}
\hspace{10pt}
a=b=1,\\
0 & \text{otherwise}
\end{cases}
\end{align*}for $a,b \in \{0,1,2\}.$
\end{lemma}\qed

\subsubsection{Computing class $\beta_{0}$, level
$(0,1)$ and $(1,0)$, and also class $\beta_{0}+f$, level $(0,-1)$
and $(-1,0)$ tubes:}We do the calculations for class $\beta_{0}$,
level $(0,1)$ and class $\beta_{0}+f$, level $(0,-1)$ tubes, the
cases where levels are on the second line bundle are similar.

We attach two tubes of levels $(0,-1)$ and $(0,1)$ to get a tube
of level $(0,0)$ (see the picture). Now applying gluing formula
and using Lemma \ref{lemma:level -1 tube} and Lemma
\ref{lemma:level 1 tube}, we get

\begin{center}
\includegraphics{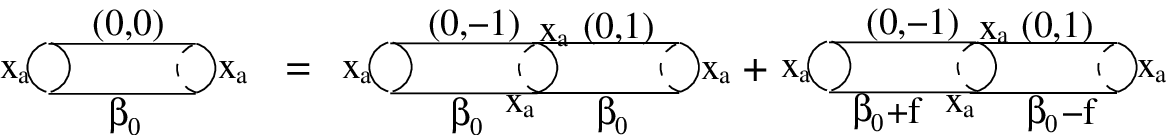}
\end{center}

\begin{align*}Z_{\beta_{0}}(0\operatorname{|}0,0)_{x_{a}x_{a}}&=
Z_{\beta_{0}}(0\operatorname{|}0,-1)_{x_{a}x_{a}}
Z_{\beta_{0}}(0\operatorname{|}0,1)_{x_{a}}^{x_{a}}\\
&\hspace{.5mm}+Z_{\beta_{0}+f}(0\operatorname{|}0,-1)_{x_{a}x_{a}}
Z_{\beta_{0}-f}(0\operatorname{|}0,1)_{x_{a}}^{x_{a}}
\end{align*} for $a=0,1,2$.\\
Again by Lemmas \ref{lemma:level -1 tube}, \ref{lemma:level 1
tube} and \ref{lemma:level 0 cap}

\begin{eqnarray*}
Z_{\beta_{0}-f}(0\operatorname{|}0,1)_{x_{0}x_{0}}&=&0,
\\Z_{\beta_{0}-f}(0\operatorname{|}0,1)_{x_{1}x_{1}}&=&0,
\\Z_{\beta_{0}}(0\operatorname{|}0,-1)_{x_{2}x_{2}}&=&0,
\end{eqnarray*}
so we can solve the equations above for the other unknowns:
\begin{eqnarray}\label{equ:1st tube reln}
Z_{\beta_{0}}(0\operatorname{|}0,1)_{x_{0}x_{0}}&=&(t_{0}-t_{1})\phi,
\nonumber
\\Z_{\beta_{0}-f}(0\operatorname{|}0,1)_{x_{1}x_{1}}&=&(t_{1}-t_{0})\phi,\nonumber
\\Z_{\beta_{0}+f}(0\operatorname{|}0,-1)_{x_{2}x_{2}}&=&\phi^{2}.
\end{eqnarray}
\\
By changing relative conditions, we can get more relations:
\begin{center}
\includegraphics{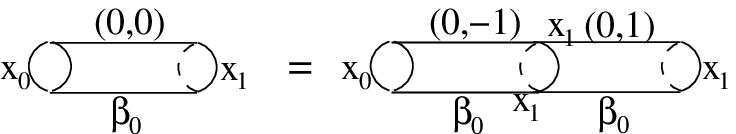}
\end{center}
\begin{align*}
0=Z_{\beta_{0}}(0\operatorname{|}0,0)_{x_{0}x_{1}}=
Z_{\beta_{0}}(0\operatorname{|}0,-1)_{x_{0}x_{1}}
Z_{\beta_{0}}(0\operatorname{|}0,1)_{x_{1}}^{x_{1}},
\end{align*}
which implies that
\begin{equation}
Z_{\beta_{0}}(0\operatorname{|}0,1)_{x_{0}x_{1}}=0.
\end{equation}
\\
We can also write
\begin{center}
\includegraphics{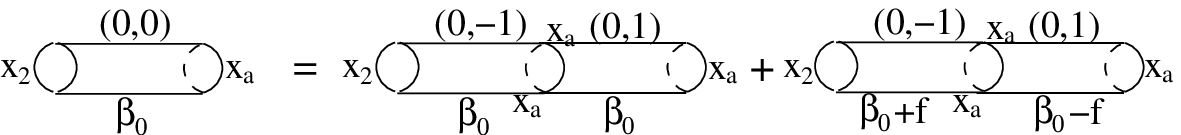}
\end{center}
\begin{align*}
0=Z_{\beta_{0}}(0\operatorname{|}0,0)_{x_{2}x_{a}}&=
Z_{\beta_{0}}(0\operatorname{|}0,-1)_{x_{2}x_{a}}
Z_{\beta_{0}}(0\operatorname{|}0,1)_{x_{a}}^{x_{a}}\\
&\hspace{.5mm}+Z_{\beta_{0}+f}(0\operatorname{|}0,-1)_{x_{2}x_{a}}
Z_{\beta_{0}-f}(0\operatorname{|}0,1)_{x_{2}}^{x_{a}}
\end{align*} for $a=0,1$. This implies that
\begin{align} \label{equ:reln tube}
&Z_{\beta_{0}}(0\operatorname{|}0,1)_{x_{0}x_{2}}=
Z_{\beta_{0}+f}(0\operatorname{|}0,-1)_{x_{0}x_{2}}(t_{2}-t_{1})\phi^{-1},\nonumber
\\&Z_{\beta_{0}}(0\operatorname{|}0,1)_{x_{1}x_{2}}=
Z_{\beta_{0}+f}(0\operatorname{|}0,-1)_{x_{1}x_{2}}(t_{2}-t_{0})\phi^{-1}.
\end{align}

Attaching the level $(0,0)$ cap to the level $(0,1)$ tube, we get
three relations:
\begin{center}
\includegraphics{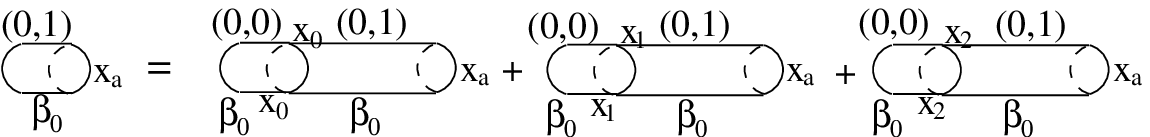}
\end{center}
\begin{align*}
0=Z_{\beta_{0}}(0\operatorname{|}0,1)_{x_{a}}^{x_{0}}+
Z_{\beta_{0}}(0\operatorname{|}0,1)_{x_{a}}^{x_{1}}+
Z_{\beta_{0}}(0\operatorname{|}0,1)_{x_{a}}^{x_{2}}
\end{align*} for $a=0,1,2$. We already know that
$$Z_{\beta_{0}}(0\operatorname{|}0,1)_{x_{0}}^{x_{1}}=
Z_{\beta_{0}}(0\operatorname{|}0,1)_{x_{1}}^{x_{0}}=0,$$ so we get
\begin{align}
&Z_{\beta_{0}}(0\operatorname{|}0,1)_{x_{0}x_{2}}=(t_{2}-t_{1})\phi,
\nonumber
\\&Z_{\beta_{0}}(0\operatorname{|}0,1)_{x_{1}x_{2}}=(t_{2}-t_{0})\phi,\nonumber
\\&Z_{\beta_{0}}(0\operatorname{|}0,1)_{x_{2}x_{2}}=(2t_{2}-t_{0}-t_{1})\phi.
\end{align}\\ Combining with (\ref{equ:reln tube})
we find
\begin{align}
&Z_{\beta_{0}+f}(0\operatorname{|}0,-1)_{x_{0}x_{2}}=\phi^{2},
\nonumber
\\&Z_{\beta_{0}+f}(0\operatorname{|}0,-1)_{x_{1}x_{2}}=\phi^{2}.
\end{align}

Finally, in order to find
$$
\begin{array} {lll}
Z_{\beta_{0}+f}(0\operatorname{|}0,-1)_{x_{0}x_{0}},
&Z_{\beta_{0}+f}(0\operatorname{|}0,-1)_{x_{1}x_{1}},
&Z_{\beta_{0}+f}(0\operatorname{|}0,-1)_{x_{0}x_{1}},
\end{array}$$
we attach the level $(0,1)$ tube to the level $(0,-1)$ tube to
obtain the class $\beta_{0}+f$, level $(0,0)$ tube:

\begin{center}
\includegraphics{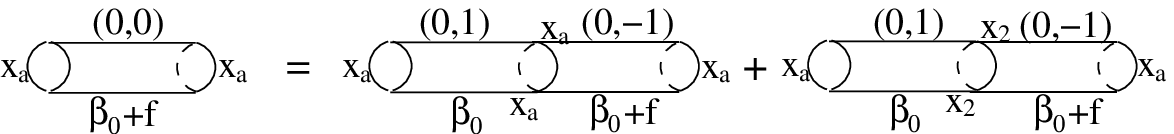}
\end{center}

\begin{align*}0=Z_{\beta_{0}+f}(0\operatorname{|}0,0)_{x_{a}x_{a}}&=
Z_{\beta_{0}}(0\operatorname{|}0,1)_{x_{a}x_{a}}
Z_{\beta_{0}+f}(0\operatorname{|}0,-1)_{x_{a}}^{x_{a}}\\
&\hspace{.5mm}+Z_{\beta_{0}}(0\operatorname{|}0,1)_{x_{a}x_{2}}
Z_{\beta_{0}+f}(0\operatorname{|}0,-1)_{x_{a}}^{x_{2}}
\end{align*} for $a=0,1$.
This implies that
\begin{align}
&Z_{\beta_{0}+f}(0\operatorname{|}0,-1)_{x_{0}x_{0}}=\phi^{2},\nonumber
\\&Z_{\beta_{0}+f}(0\operatorname{|}0,-1)_{x_{1}x_{1}}=\phi^{2},
\end{align}
And similarly,
\begin{center}
\includegraphics{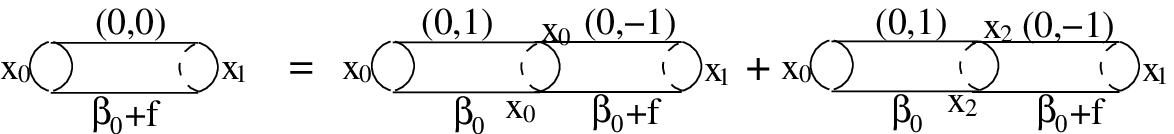}
\end{center}

\begin{align*}0=Z_{\beta_{0}+f}(0\operatorname{|}0,0)_{x_{0}x_{1}}&=
Z_{\beta_{0}}(0\operatorname{|}0,1)_{x_{0}x_{0}}
Z_{\beta_{0}+f}(0\operatorname{|}0,-1)_{x_{0}}^{x_{1}}\\
&\hspace{.5mm}+Z_{\beta_{0}}(0\operatorname{|}0,1)_{x_{0}x_{2}}
Z_{\beta_{0}+f}(0\operatorname{|}0,-1)_{x_{1}}^{x_{2}},
\end{align*}
which implies that
\begin{align} \label{equ:last tube reln}
Z_{\beta_{0}+f}(0\operatorname{|}0,-1)_{x_{0}x_{1}}=\phi^{2}.
\end{align}

We now summarize (\ref{equ:1st tube reln}-\ref{equ:last tube
reln}) into the following lemma:
\begin{lemma} \label{lemma:level 1 , -1 tube}
Partition functions for the class $\beta_{0}$, level $(0,1)$ and
$(1,0)$ tubes, and also for the class $\beta_{0}+f$, level
$(0,-1)$ and $(-1,0)$ tubes are given by
\begin{align*}
&\left[Z_{\beta_{0}}(0\operatorname{|}0,1)_{x_{a}x_{b}}\right]\text{\hspace{.6cm}}=
\left[
\begin{array}{ccc}
t_{0}-t_{1} & 0 & t_{2}-t_{1} \\
0 & t_{1}-t_{0} & t_{2}-t_{0} \\
t_{2}-t_{1} & t_{2}-t_{0} & 2t_{2}-t_{0}-t_{1}
\end{array} \right]\phi,\\
&\left[Z_{\beta_{0}}(0\operatorname{|}1,0)_{x_{a}x_{b}}\right]\text{\hspace{.6cm}}=
\left[
\begin{array}{ccc}
t_{0}-t_{2} & t_{1}-t_{2} & 0 \\
t_{1}-t_{2} & 2t_{1}-t_{0}-t_{2} & t_{1}-t_{0} \\
0 & t_{1}-t_{0} & t_{2}-t_{0}
\end{array} \right]\phi,\\
&\left[Z_{\beta_{0}+f}(0\operatorname{|}0,-1)_{x_{a}x_{b}}\right]=
\left[\begin{array}{ccc}
1 & 1 & 1 \\
1 & 1 & 1 \\
1 & 1 & 1
\end{array} \right]\phi^{2},\\
&\left[Z_{\beta_{0}+f}(0\operatorname{|}-1,0)_{x_{a}x_{b}}\right]=\left[\begin{array}{lll}
1 & 1 & 1 \\
1 & 1 & 1 \\
1 & 1 & 1
\end{array} \right]\phi^{2},
\end{align*}
for $a,b\in \{0,1,2\}$, where partition functions with the index
$x_{a}x_{b}$ are the $(a+1,b+1)$ entry of the matrices above.
\end{lemma}\qed
\subsubsection{computing class $\beta_{0}+f$, level
$(0,0)$ pants:}
\begin{lemma} \label{lemma:level 0 pants}
Partition functions for the class $\beta_{0}+f$, level $(0,0)$
pants are given by
\begin{align*}
&Z_{\beta_{0}+f}(0\operatorname{|}0,0)_{x_{0}x_{1}x_{2}}=0,\\
&Z_{\beta_{0}+f}(0\operatorname{|}0,0)_{x_{0}x_{2}x_{2}}=(t_{2}-t_{1})\phi^{3},\\
&Z_{\beta_{0}+f}(0\operatorname{|}0,0)_{x_{1}x_{2}x_{2}}=(t_{2}-t_{0})\phi^{3},\\
&Z_{\beta_{0}+f}(0\operatorname{|}0,0)_{x_{0}x_{0}x_{2}}=(t_{1}-t_{0})\phi^{3},\\
&Z_{\beta_{0}+f}(0\operatorname{|}0,0)_{x_{1}x_{1}x_{2}}=(t_{0}-t_{1})\phi^{3},\\
&Z_{\beta_{0}+f}(0\operatorname{|}0,0)_{x_{2}x_{2}x_{2}}=(2t_{2}-t_{0}-t_{1})\phi^{3},\\
&Z_{\beta_{0}+f}(0\operatorname{|}0,0)_{x_{0}x_{1}x_{1}}=(t_{1}-t_{2})\phi^{3},\\
&Z_{\beta_{0}+f}(0\operatorname{|}0,0)_{x_{0}x_{0}x_{1}}=(t_{0}-t_{2})\phi^{3},\\
&Z_{\beta_{0}+f}(0\operatorname{|}0,0)_{x_{0}x_{0}x_{0}}=(2t_{0}-t_{1}-t_{2})\phi^{3},\\
&Z_{\beta_{0}+f}(0\operatorname{|}0,0)_{x_{1}x_{1}x_{1}}=(2t_{1}-t_{0}-t_{2})\phi^{3},
\end{align*}
for $a,b\in \{0,1,2\}$.
\end{lemma}
\textsc{Proof:} We attach the level $(0,1)$ cap to the level
$(0,0)$ pants to obtain the class $\beta_{0}$, level $(0,1)$ tube.
Applying the gluing formula together with Lemma \ref{lemma:level 1
cap} and \ref{lemma:level 1 , -1 tube}, we get the following
relation
\begin{center}
\includegraphics{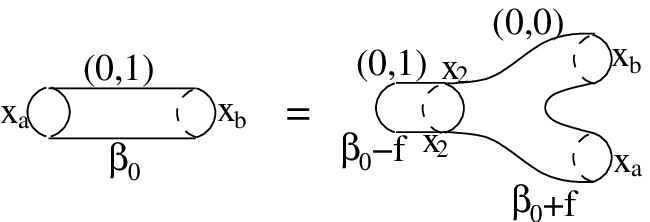}
\end{center}

$$Z_{\beta_{0}+f}(0\operatorname{|}0,0)_{x_{2}x_{a}x_{b}}=
Z_{\beta_{0}+f}(0\operatorname{|}0,1)_{x_{a}x_{b}}\phi^{2}.$$ From
this we can get all
$Z_{\beta_{0}+f}(0\operatorname{|}0,0)_{x_{a}x_{b}x_{c}}$ with at
least one of $a,b,c$ is equal to $2$.

If we attach the level $(0,-1)$ cap to the level $(0,0)$ pants to
obtain the class $\beta+f$, level $(0,-1)$ tube we will get

\begin{center}
\includegraphics{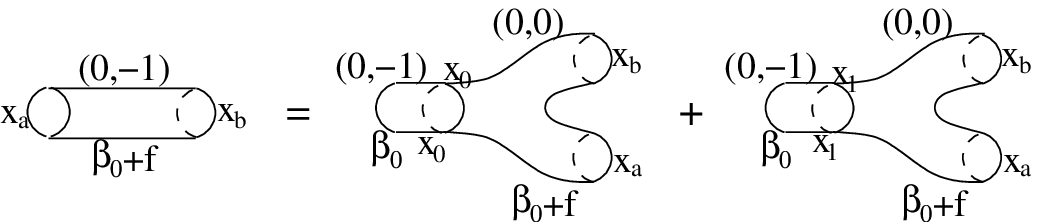}
\end{center}

\begin{align} \label{equ:pants reln}
&Z_{\beta_{0}+f}(0\operatorname{|}0,0)_{x_{0}x_{0}x_{1}}-
Z_{\beta_{0}+f}(0\operatorname{|}0,0)_{x_{0}x_{1}x_{1}}
=(t_{0}-t_{1})\phi^{3},\nonumber
\\&Z_{\beta_{0}+f}(0\operatorname{|}0,0)_{x_{0}x_{0}x_{0}}-
Z_{\beta_{0}+f}(0\operatorname{|}0,0)_{x_{0}x_{0}x_{1}}
=(t_{0}-t_{1})\phi^{3},\nonumber
\\&Z_{\beta_{0}+f}(0\operatorname{|}0,0)_{x_{0}x_{1}x_{1}}-
Z_{\beta_{0}+f}(0\operatorname{|}0,0)_{x_{1}x_{1}x_{1}}
=(t_{0}-t_{1})\phi^{3}.
\end{align}
We now write a Frobenius relation as follows
\begin{center}
\includegraphics{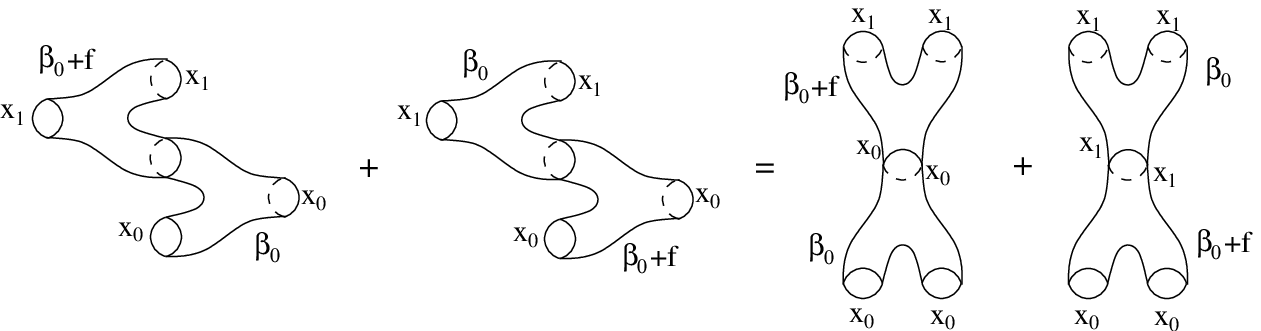}
\end{center}
$$0=
Z_{\beta_{0}+f}(0\operatorname{|}0,0)_{x_{0}x_{1}x_{1}}
Z_{\beta_{0}}(0\operatorname{|}0,0)_{x_{0}x_{0}}^{x_{0}}+
Z_{\beta_{0}}(0\operatorname{|}0,0)_{x_{1}x_{1}x_{1}}
Z_{\beta_{0}+f}(0\operatorname{|}0,0)_{x_{0}x_{0}}^{x_{1}},$$
where the left hand side is zero by Lemma {\ref{lemma:level 0
cap}.\\By Lemma \ref{lemma:level 0 cap}, this simplifies to
$$(t_{0}-t_{2})Z_{\beta_{0}+f}(0\operatorname{|}0,0)_{x_{0}x_{1}x_{1}}=
(t_{1}-t_{2})Z_{\beta_{0}+f}(0\operatorname{|}0,0)_{x_{0}x_{0}x_{1}}.$$
Combining this with (\ref{equ:pants reln}), we will find the rest
of the partition functions in the lemma.\qed

We now know all the partition functions of pants, so we are able
to prove the semisimplicity of the TQFT (see Remark
\ref{rem:TQFT}):

\begin{prop} \label{pro:semisimple}
The level $(0,0)$ TQFT, resulted from our setting and Theorem
\ref{thm:gluing rules} is semisimple.
\end{prop}
\textsc{Proof:} By Lemma \ref{lemma:level 0 cap} and
\ref{lemma:level 0 pants}, we have
$$Z(0\operatorname{|}0,0)_{x_{a}x_{b}}^{x_{c}}|_{u=0}=
Z_{\beta_{0}}(0\operatorname{|}0,0)_{x_{a}x_{b}}^{x_{c}}=
\begin{cases} T(x_{a}) & \text{if}
\hspace{10pt} a=b=c, \\ 0 & \text{otherwise.}
\end{cases}
$$
This means that for $u=0$, the basis
$\{\frac{e_{x_{0}}}{T(x_{0})},\frac{e_{x_{1}}}{T(x_{1})},\frac{e_{x_{2}}}{T(x_{2})}\}$
of the corresponding Frobenius algebra (see Remark \ref{rem:TQFT})
is idempotent:
$$\frac{e_{x_{i}}}{T(x_{i})}\otimes \frac{e_{x_{j}}}{T(x_{j})}=
\delta_{i}^{j}\frac{e_{x_{i}}}{T(x_{i})}$$
 This proves the semisimplicity when $u=0$ (see
\cite{Br-Pa-TQFT}, Section 2). Now the proposition follows from
Proposition 2.2 in \cite{Br-Pa-TQFT}. \qed

\section {Proof of Theorem \ref{thm:matrix formula}} \label{sec:result}
We now know everything we need in order to prove Theorem
\ref{thm:matrix formula}. We first find the \textbf{first} and
\textbf{second} \textbf{level creation operators} and also
\textbf{genus adding operator} which are by definition
\begin{align*}
&U_{1}=\left[Z(0\operatorname{|}1,0)_{x_{a}}^{x_{b}}\right],\\
&U_{2}=\left[Z(0\operatorname{|}0,1)_{x_{a}}^{x_{b}}\right],\\
&G\text{\hspace{1.3mm}}=\left[Z(1\operatorname{|}0,0)_{x_{a}}^{x_{b}}\right],
\end{align*} respectively. Here partition functions with the lower
index $x_{a}$ and the upper index $x_{b}$ are the $(a+1,b+1)$
entry of the matrices above.

We can find $U_{1}$ and $U_{2}$ by simply raising the indices in
Lemma  \ref{lemma:level 1 tube} and \ref{lemma:level 1 , -1 tube}:
\begin{align*}
U_{1}&=\left[Z_{\beta_{0}-f}(0\operatorname{|}1,0)_{x_{a}}^{x_{b}}\right]+
\left[Z_{\beta_{0}}(0\operatorname{|}1,0)_{x_{a}}^{x_{b}}\right]\\
&=\left[
\begin{array}{ccc}
0 & 0 & 0 \\
0 &
(t_{1}-t_{0})(t_{1}-t_{2}) & 0 \\
0 & 0 & 0
\end{array} \right]\phi^{-2}+\left[
\begin{array}{ccc}
\frac{1}{t_{0}-t_{1}} &
\frac{t_{1}-t_{2}}{(t_{0}-t_{1})(t_{0}-t_{2})}
& 0 \\
\frac{1}{t_{1}-t_{0}} &
\frac{2t_{1}-t_{0}-t_{2}}{(t_{1}-t_{0})(t_{1}-t_{2})} &
\frac{1}{t_{1}-t_{2}} \\
0 & \frac{t_{1}-t_{0}}{(t_{2}-t_{0})(t_{2}-t_{1})} &
\frac{1}{t_{2}-t_{1}}
\end{array} \right]\phi\\
&=\left[
\begin{array}{ccc}
\frac{\phi}{t_{0}-t_{1}} &
\frac{(t_{1}-t_{2})\phi}{(t_{0}-t_{1})(t_{0}-t_{2})} & 0 \\
\frac{\phi}{t_{1}-t_{0}} &
\frac{(t_{1}-t_{0})^{2}(t_{1}-t_{2})^{2}\phi^{-2}+(2t_{1}-t_{0}-t_{2})\phi}
{(t_{1}-t_{0})(t_{1}-t_{2})} & \frac{\phi}{t_{1}-t_{2}} \\
0 & \frac{(t_{1}-t_{0})\phi}{(t_{2}-t_{0})(t_{2}-t_{1})} &
\frac{\phi}{t_{2}-t_{1}}
\end{array} \right],
\end{align*}
\begin{align} \label{equ:level creation operators}
U_{2}&=\left[Z_{\beta_{0}-f}(0\operatorname{|}0,1)_{x_{a}}^{x_{b}}\right]+
\left[Z_{\beta_{0}}(0\operatorname{|}0,1)_{x_{a}}^{x_{b}}\right]\nonumber\\
&=\left[
\begin{array}{ccc}
0 & 0 & 0 \\
0 &
0 & 0 \\
0 & 0 & (t_{2}-t_{0})(t_{2}-t_{1})
\end{array} \right]\phi^{-2}+\left[
\begin{array}{ccc}
\frac{1}{t_{0}-t_{2}} & 0 &
\frac{t_{2}-t_{1}}{(t_{0}-t_{1})(t_{0}-t_{2})} \\
0 & \frac{1}{t_{1}-t_{2}} &
\frac{t_{2}-t_{0}}{(t_{1}-t_{0})(t_{1}-t_{2})} \\
\frac{1}{t_{2}-t_{0}} & \frac{1}{t_{2}-t_{1}} &
\frac{2t_{2}-t_{0}-t_{1}}{(t_{2}-t_{0})(t_{2}-t_{1})}
\end{array} \right]\phi \nonumber\\
&=\left[
\begin{array}{ccc}
\frac{\phi}{t_{0}-t_{2}} & 0 &
\frac{(t_{2}-t_{1})\phi}{(t_{0}-t_{1})(t_{0}-t_{2})} \\
0 & \frac{\phi}{t_{1}-t_{2}} &
\frac{(t_{2}-t_{0})\phi}{(t_{1}-t_{0})(t_{1}-t_{2})} \\
\frac{\phi}{t_{2}-t_{0}} & \frac{\phi}{t_{2}-t_{1}} &
\frac{(t_{2}-t_{0})^{2}(t_{2}-t_{1})^{2}\phi^{-2}+
(2t_{2}-t_{0}-t_{1})\phi}{(t_{2}-t_{0})(t_{2}-t_{1})}
\end{array} \right].
\end{align}
By Lemma \ref{lemma:level 0 cap},
$\left[Z_{\beta_{0}}(0\operatorname{|}0,0)_{x_{a}}^{x_{b}}\right]$
is the identity matrix, and by the gluing formula we can write
\begin{align*}
&\left[Z(0\operatorname{|}0,0)_{x_{a}}^{x_{b}}\right]=
\left[Z(0\operatorname{|}0,1)_{x_{a}}^{x_{b}}\right]
\left[Z(0\operatorname{|}0,-1)_{x_{a}}^{x_{b}}\right],\\
&\left[Z(0\operatorname{|}0,0)_{x_{a}}^{x_{b}}\right]=
\left[Z(0\operatorname{|}1,0)_{x_{a}}^{x_{b}}\right]
\left[Z(0\operatorname{|}-1,0)_{x_{a}}^{x_{b}}\right].
\end{align*}
Therefore, by Lemma \ref{lemma:level -1 tube} and \ref{lemma:level
1 , -1 tube} we have
\begin{align*}
U_{1}^{-1}&=\left[Z_{\beta_{0}}(0\operatorname{|}-1,0)_{x_{a}}^{x_{b}}\right]+
\left[Z_{\beta_{0}+f}(0\operatorname{|}-1,0)_{x_{a}}^{x_{b}}\right]\\
&=\left[\begin{array}{ccc}
t_{0}-t_{1} & 0 & 0 \\
0 &
0 & 0 \\
0 & 0 & t_{2}-t_{1}
\end{array} \right]\phi^{-1}\\
&+\left[\begin{array}{ccc} \frac{1}{(t_{0}-t_{1})(t_{0}-t_{2})} &
\frac{1}{(t_{0}-t_{1})(t_{0}-t_{2})}
& \frac{1}{(t_{0}-t_{1})(t_{0}-t_{2})} \\
\frac{1}{(t_{1}-t_{0})(t_{1}-t_{2})} &
\frac{1}{(t_{1}-t_{0})(t_{1}-t_{2})}
& \frac{1}{(t_{1}-t_{0})(t_{1}-t_{2})} \\
\frac{1}{(t_{2}-t_{0})(t_{2}-t_{1})} &
\frac{1}{(t_{2}-t_{0})(t_{2}-t_{1})} &
\frac{1}{(t_{2}-t_{0})(t_{2}-t_{1})}
\end{array} \right]\phi^{2}\\
&=\left[\begin{array}{ccc}
\frac{(t_{0}-t_{1})^{2}(t_{0}-t_{2})\phi^{-1}+\phi^{2}}{(t_{0}-t_{1})(t_{0}-t_{2})}
& \frac{\phi^{2}}{(t_{0}-t_{1})(t_{0}-t_{2})}
& \frac{\phi^{2}}{(t_{0}-t_{1})(t_{0}-t_{2})} \\
\frac{\phi^{2}}{(t_{1}-t_{0})(t_{1}-t_{2})} &
\frac{\phi^{2}}{(t_{1}-t_{0})(t_{1}-t_{2})}
& \frac{\phi^{2}}{(t_{1}-t_{0})(t_{1}-t_{2})} \\
\frac{\phi^{2}}{(t_{2}-t_{0})(t_{2}-t_{1})} &
\frac{\phi^{2}}{(t_{2}-t_{0})(t_{2}-t_{1})} &
\frac{(t_{2}-t_{0})(t_{2}-t_{1})^{2}\phi^{-1}+\phi^{2}}{(t_{2}-t_{0})(t_{2}-t_{1})}
\end{array} \right],
\end{align*}
\begin{align} \label{equ:level anlhn operators}
U_{2}^{-1}&=\left[Z_{\beta_{0}}(0\operatorname{|}0,-1)_{x_{a}}^{x_{b}}\right]+
\left[Z_{\beta_{0}+f}(0\operatorname{|}0,-1)_{x_{a}}^{x_{b}}\right]\nonumber\\
&=\left[
\begin{array}{ccc}
t_{0}-t_{2} & 0 & 0 \\
0 &
t_{1}-t_{2} & 0 \\
0 & 0 & 0
\end{array}\right]\phi^{-1} \nonumber \\
&+\left[\begin{array}{ccc} \frac{1}{(t_{0}-t_{1})(t_{0}-t_{2})} &
\frac{1}{(t_{0}-t_{1})(t_{0}-t_{2})}
& \frac{1}{(t_{0}-t_{1})(t_{0}-t_{2})} \\
\frac{1}{(t_{1}-t_{0})(t_{1}-t_{2})} &
\frac{1}{(t_{1}-t_{0})(t_{1}-t_{2})}
& \frac{1}{(t_{1}-t_{0})(t_{1}-t_{2})} \\
\frac{1}{(t_{2}-t_{0})(t_{2}-t_{1})} &
\frac{1}{(t_{2}-t_{0})(t_{2}-t_{1})} &
\frac{1}{(t_{2}-t_{0})(t_{2}-t_{1})}
\end{array} \right]\phi^{2} \nonumber\\
&=\left[\begin{array}{ccc}
\frac{(t_{0}-t_{1})(t_{0}-t_{2})^{2}\phi^{-1}+\phi^{2}}{(t_{0}-t_{1})(t_{0}-t_{2})}
& \frac{\phi^{2}}{(t_{0}-t_{1})(t_{0}-t_{2})}
& \frac{\phi^{2}}{(t_{0}-t_{1})(t_{0}-t_{2})} \\
\frac{\phi^{2}}{(t_{1}-t_{0})(t_{1}-t_{2})} &
\frac{(t_{1}-t_{0})(t_{1}-t_{2})^{2}\phi^{-1}+\phi^{2}}{(t_{1}-t_{0})(t_{1}-t_{2})}
& \frac{\phi^{2}}{(t_{1}-t_{0})(t_{1}-t_{2})} \\
\frac{\phi^{2}}{(t_{2}-t_{0})(t_{2}-t_{1})} &
\frac{\phi^{2}}{(t_{2}-t_{0})(t_{2}-t_{1})} &
\frac{\phi^{2}}{(t_{2}-t_{0})(t_{2}-t_{1})}
\end{array} \right].
\end{align}
$U_{1}^{-1}$ and $U_{2}^{-1}$ are the \textbf{first} and
\textbf{second level annihilation operators}, respectively.

Now we are going to find the matrix $G$. By the same argument as
one given for Lemma \ref{lem:tqft pieces} one can prove that
$$Z(1\operatorname{|}0,0)_{x_{a}x_{b}}=Z_{\beta_{0}}(1\operatorname{|}0,0)_{x_{a}x_{b}}+
Z_{\beta_{0}+f}(1\operatorname{|}0,0)_{x_{a}x_{b}}$$ for $a,b \in
\{0,1,2\}$. Thus we have
$$G=\left[Z_{\beta_{0}}(1\operatorname{|}0,0)_{x_{a}}^{x_{b}}\right]+
\left[Z_{\beta_{0}+f}(1\operatorname{|}0,0)_{x_{a}}^{x_{b}}\right].$$
For calculating the terms in the right hand side of this, we
attach two pants at two points (see the picture), and apply the
gluing formula:
\begin{center}
\includegraphics{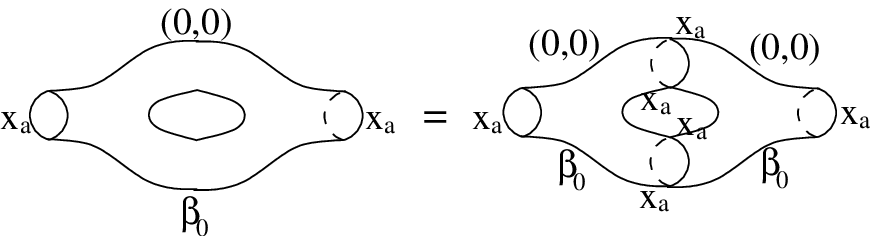}
\end{center}
\begin{center}
\includegraphics{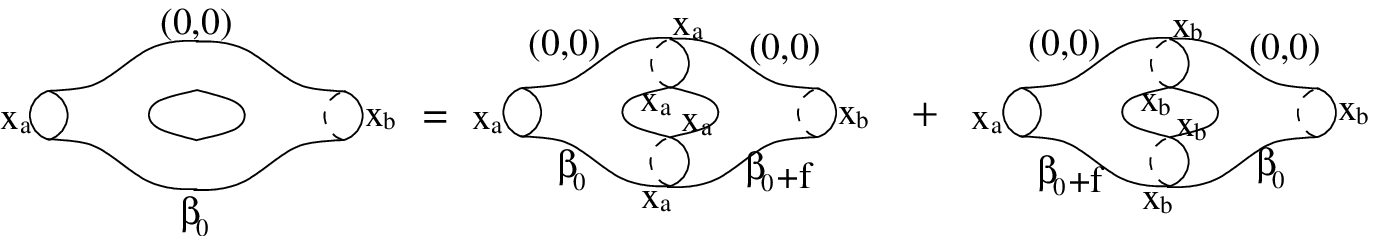}
\end{center}
\begin{align*}
&Z_{\beta_{0}}(1\operatorname{|}0,0)_{x_{a}x_{a}}\text{\hspace{.35cm}}=
Z_{\beta_{0}}(0\operatorname{|}0,0)_{x_{a}x_{a}x_{a}}
Z_{\beta_{0}}(0\operatorname{|}0,0)_{x_{a}}^{x_{a}x_{a}},\\
&Z_{\beta_{0}+f}(1\operatorname{|}0,0)_{x_{a}x_{b}}=
Z_{\beta_{0}}(0\operatorname{|}0,0)_{x_{a}x_{a}x_{a}}
Z_{\beta_{0}+f}(0\operatorname{|}0,0)_{x_{b}}^{x_{a}x_{a}}\\&
\text{\hspace{3.1cm}}+Z_{\beta_{0}+f}(0\operatorname{|}0,0)_{x_{a}x_{b}x_{b}}
Z_{\beta_{0}}(0\operatorname{|}0,0)_{x_{b}}^{x_{b}x_{b}},
\end{align*} which implies that

\begin{equation} \label{equ:beta0 g adding}
\left[Z_{\beta_{0}}(1\operatorname{|}0,0)_{x_{a}}^{x_{b}}\right]=
\left[\begin{array}{ccc} (t_{0}-t_{1})(t_{0}-t_{2}) & 0
& 0 \\
0 & (t_{1}-t_{0})(t_{1}-t_{2})
& 0 \\
0 & 0 & (t_{2}-t_{0})(t_{2}-t_{1}))
\end{array}
\right]
\end{equation} and
\begin{equation} \label{equ:beta0+f g adding}
\left[Z_{\beta_{0}+f}(1\operatorname{|}0,0)_{x_{a}}^{x_{b}}\right]=
\left[\begin{array}{ccc}
\frac{2(2t_{0}-t_{1}-t_{2})}{(t_{0}-t_{1})(t_{0}-t_{2})} &
\frac{t_{0}+t_{1}-2t_{2}}{(t_{0}-t_{1})(t_{0}-t_{2})}
&\frac{t_{0}+t_{2}-2t_{1}}{(t_{0}-t_{1})(t_{0}-t_{2})} \\
\frac{t_{0}+t_{1}-2t_{2}}{(t_{1}-t_{0})(t_{1}-t_{2})} &
\frac{2(2t_{1}-t_{0}-t_{2})}{(t_{1}-t_{0})(t_{1}-t_{2})}
&\frac{t_{1}+t_{2}-2t_{0}}{(t_{1}-t_{0})(t_{1}-t_{2})} \\
\frac{t_{0}+t_{2}-2t_{1}}{(t_{2}-t_{0})(t_{2}-t_{1})} &
\frac{t_{1}+t_{2}-2t_{0}}{(t_{2}-t_{0})(t_{2}-t_{1})} &
\frac{2(2t_{2}+t_{0}-t_{1})}{(t_{2}-t_{0})(t_{2}-t_{1})}
\end{array} \right]\phi^{3}.
\end{equation}

Now to prove the formula in Theorem \ref{thm:matrix formula}, we
glue a chain of $g-1$ genus one, level $(0,0)$ cobordisms from a
circle to another circle, which are attached subsequently at their
ends to a chain of $k_{1}$ level $(1,0)$ and $k_{2}$ level $(0,1)$
tubes, and after all, two ends of the resulting chain are also
attached to each other (see the picture below).
\begin{center}
\includegraphics{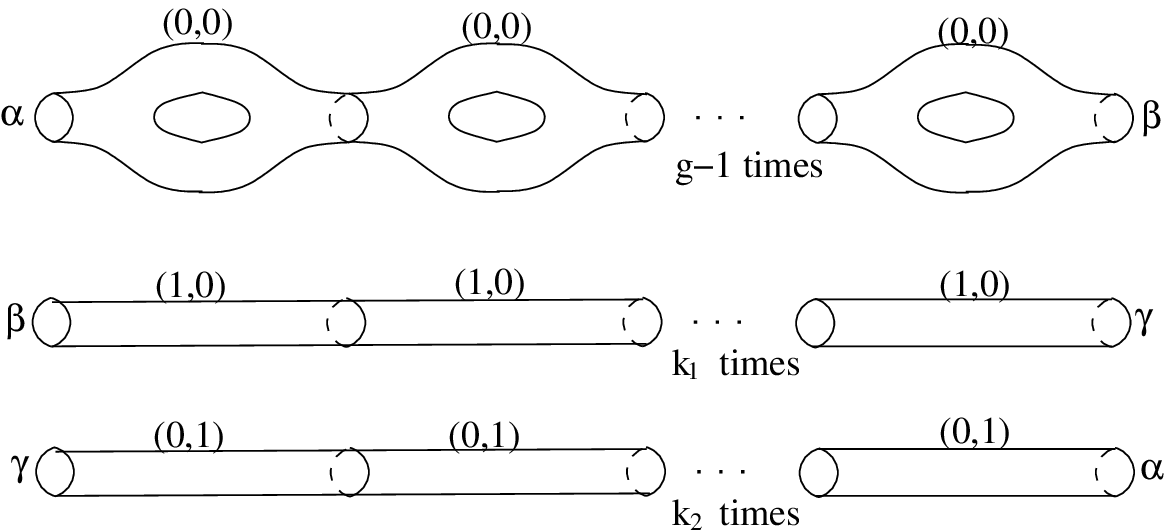}
\end{center}
In this picture, objects with the same mark at their end are
attached to each other. Applying the gluing formula each time we
glue along a fiber in the procedure above, we obtain the formula
in Theorem \ref{thm:matrix formula}:
$$Z(g\operatorname{|}k_{1},k_{2})=
\operatorname{tr}\left(G^{g-1}U_{1}^{k_{1}}U_{2}^{k_{2}}\right).$$
Thus we have proven Theorem \ref{thm:matrix formula} in the case
where $g\geq 1$. The same formula holds for $g=0$, which follows
from the semisimplicity of the level $(0,0)$ TQFT (see Remark
\ref{rem:TQFT}, and Theorem \ref{pro:semisimple}).

\section {Proof of Theorem \ref{thm:Calabi-Yau class}} \label{sec:proof
of theorem} We begin by outlining our proof of Theorem
\ref{thm:Calabi-Yau class}:

We first prove that any $\mathbb{P}^{2}$-bundle over a curve $C$
is deformation equivalent to
$$\mathbb{P}(\mathcal{O}\oplus\mathcal{O}\oplus L)\rightarrow C.$$
This is the space for which we obtained a formula for finding the
section class equivariant Gromov-Witten invariants (Theorem
\ref{thm:matrix formula}). We will expand this formula and look
for those terms which are corresponding to the Calabi-Yau section
classes. There are lots of such terms, but we will see that all of
them except three terms are zero for some reason; all equivariant
variable will cancel each other and these three terms will add up
to give the formula in the theorem.

We first prove the following lemma:
\begin{lem} \label{lem:deformation}
Any $\mathbb{P}^2$-bundle over a curve $C$ is deformation
equivalent to
$$\mathbb{P}(\mathcal{O}\oplus\mathcal{O}\oplus L)\rightarrow C,$$
where $\mathcal{O}\rightarrow C$ is the trivial bundle and
$L\rightarrow C$ is a line bundle of degree $k$.
\end{lem}
\textsc{Proof:} First, we show that every $\mathbb{P}^2$-bundle
over a curve $C$ is of the form $$\mathbb{P}(E)\rightarrow C,$$
where $E\rightarrow C$ is a rank $3$ bundle.

A rank three vector bundle (resp. $\mathbb{P}^2$-bundle) over $C$
is classified by an element in $\check{H}^{1}(C,Gl(3))$ (resp.
$\check{H}^{1}(C,PGl(3))$), where $Gl(3)$ (resp. $PGl(3)$) is the
(non-abelian) sheaf of $Gl(3)$ (resp. $PGl(3)$) valued holomorphic
functions on $C$. From the exact sequence of sheaves
$$0\rightarrow \mathcal{O}^{*}\rightarrow Gl(3) \rightarrow
PGl(3)\rightarrow 0,$$ we get a map
$$\check{H}^{1}(C,PGl(3))\rightarrow \check{H}^{2}(C,\mathcal{O}^{*}).$$
By examining the cocycles, one can see that a
$\mathbb{P}^2$-bundle over $C$ is of the form
$\mathbb{P}(E)\rightarrow C$ if and only if the corresponding
element in $\check{H}^{1}(C,PGl(3))$ goes to zero under the above
map. This element is represented by the $\check{\text{C}}$ech
cocycle obtained from the transition functions of the bundle. But
$$\check{H}^{2}(C,\mathcal{O}^{*})=0$$ for a curve $C$; this completes the
first part of the proof of the lemma.

Next, we show that $\mathbb{P}(E)\rightarrow C$ is deformation
equivalent to $$\mathbb{P}(\mathcal{O}\oplus\mathcal{O}\oplus
L)\rightarrow C.$$ It is an standard fact that for a rank $3$
bundle $E$ over a curve we have the following exact sequence of
bundles over $C$ (see \cite{Hu-Le}, Example $5.0.1$):
$$0\rightarrow \mathcal{O}\oplus\mathcal{O}\rightarrow E(m)
\rightarrow L\rightarrow 0,$$ for some line bundle $L$ and some
$m\gg0$ such that $E(m)$ is globally generated by its sections.
Thus $E(m)$ corresponds to an element $$v \in
\operatorname{Ext}^{1}(L,\mathcal{O}\oplus\mathcal{O}).$$ We can
deform $E(m)$ by deforming the extension class $v$ to $0$ inside
this vector space. But $0
\in\operatorname{Ext}^{1}(L,\mathcal{O}\oplus\mathcal{O})$
corresponds to
$$\mathcal{O}\oplus\mathcal{O}\oplus L\rightarrow C.$$ So we have
proved that $E(m)$ is deformation equivalent to
$\mathcal{O}\oplus\mathcal{O}\oplus L$. Now we use the isomorphism
$$\mathbb{P}(E) \cong \mathbb{P}(E(m))$$ to complete the proof of
the lemma. \qed

By this lemma, we can assume that the space $X$ in the theorem is
of the form  $$\mathbb{P}(\mathcal{O}\oplus\mathcal{O}\oplus
L)\rightarrow C.$$ For simplicity, we use the following notations
in this section:
\begin{align*}
A&=\left[Z_{\beta_{0}}(1\operatorname{|}0,0)_{x_{a}}^{x_{b}}\right]
\text{\hspace{1.8cm}given by (\ref{equ:beta0 g adding})},\\
B&=\left[Z_{\beta_{0}+f}(1\operatorname{|}0,0)_{x_{a}}^{x_{b}}\right]
\text{\hspace{1.4cm}given by (\ref{equ:beta0+f g adding})},\\
C&=\left[Z_{\beta_{0}-f}(0\operatorname{|}0,1)_{x_{a}}^{x_{b}}\right]
\text{\hspace{1.4cm}given by (\ref{equ:level creation operators})},\\
E&=\left[Z_{\beta_{0}}(0\operatorname{|}0,1)_{x_{a}}^{x_{b}}\right]
\text{\hspace{1.8cm}given by (\ref{equ:level creation operators})},\\
N&=\left[Z_{\beta_{0}}(0\operatorname{|}0,-1)_{x_{a}}^{x_{b}}\right]
\text{\hspace{1.5cm}given by (\ref{equ:level anlhn operators})},\\
M&=\left[Z_{\beta_{0}+f}(0\operatorname{|}0,-1)_{x_{a}}^{x_{b}}\right]
\text{\hspace{1.1cm}given by (\ref{equ:level anlhn operators})}.\\
\end{align*}
Then we have
\begin{align*}
&G\hspace{3.7mm}=A+B,\\ &U_{2}\hspace{2.5mm}=C+E,\\
&U_{2}^{-1}=N+M,
\end{align*}
and we can write the formula in Theorem \ref{thm:matrix formula}
for $k_{1}=0$ and $k_{2}=k$ as follows:
$$Z(g\operatorname{|}0,k)=\operatorname{tr}\left((A+B)^{g-1}(C+E)^{k}\right).$$
Now we are looking for those terms in this formula that correspond
to Calabi-Yau section class. If we denote this class by
$$\beta_{cs}=\beta_{0}+nf,$$
then $n$ must satisfy this equation:
\begin{align*}
K_{X}\cdot \beta &=0\quad \Rightarrow \\
2g-2-k-3n&=0.
\end{align*}
If for given $g$ and $k$, there is a integral solution for $n$ in
this equation then the Calabi-Yau class exists. We write the above
equation in terms of $n$ instead of $k$:
$$
Z(g\operatorname{|}0,k)=\operatorname{tr}\left((A+B)^{g-1}(C+E)^{2g-2-3n}\right).$$
Now by the gluing formula, $G=A+B$ commutes with $U_{2}=C+E$, so
we have
\begin{equation} \label{equ:z(g|0,k)}
Z(g\operatorname{|}0,k)=
\operatorname{tr}\left(\left((A+B)(C+E)^{2}\right)^{g-1}(C+E)^{-3n}\right).
\end{equation}
\textbf{Notation.} For two matrices $U$ and $V$, by
$(U^{a},V^{b})$ for $a,b \in \mathbb{Z}^{+}$, we mean the sum of
the all products that we can write containing $a$ copies of $U$
and $b$ copies of $V$. For example
$$(U^{2},V)=U^{2}V+UVU+VU^{2}.$$

We first assume that $g>0$. We distinguish two cases:
\begin{enumerate}
\item[(i)] $n<0$\\
One can see that $E^{3}=0$ and $BE^{2}=0$, so we have
\begin{align*}
&Z(g\operatorname{|}0,k)=\\&\operatorname{tr}\left(\left((A+B)\left(C^{2}+E^{2}
+(E,C)\right)\right)^{g-1}\left((E^{2},C)+(E,C^{2})+C^{3}\right)^{-n}\right)\\
&\hspace{1.65cm}=\operatorname{tr}\left(\left(AE^{2}+B(E,C)+\dots\right)^{g-1}
\left((E^{2},C)+\dots\right)^{-n}\right).
\end{align*}
$A,B,C$ and $E$ correspond to the class
$\beta_{0},\beta_{0}+f,\beta_{0}-f$, and $\beta_{0}$,
respectively. One can see that the only those terms that have been
written in the last equality above contribute to make the class
$\beta_{cs}=\beta_{0}+nf$. Thus
\begin{align}\label{equ:trace equation}
&Z_{\beta_{cs}}(g\operatorname{|}0,k)\nonumber
\\&=\operatorname{tr}\left(\left(AE^{2}+B(E,C)\right)^{g-1}
\left((E^{2},C)\right)^{-n}\right)\nonumber \\
&=\operatorname{tr}\left(\left(AE^{2}+BEC+BCE\right)^{g-1}
\left(E^{2}C+CE^{2}+ECE\right)^{-n}\right)\nonumber \\
&=\operatorname{tr}\left((AE^{2})^{g-1}(CE^{2})^{-n}\right)
+\operatorname{tr}\left((BEC)^{g-1}(E^{2}C)^{-n}\right)\nonumber \\
&\hspace{.5mm}+\operatorname{tr}\left((BCE)^{g-1}(ECE)^{-n}\right).
\end{align}
For the last equality we only used the fact that
$$E^{3}=BE^{2}=0$$ again and also
$$\operatorname{tr}(UV)=\operatorname{tr}(VU),$$ for any two
matrices $U$ and $V$.

Now one can see easily by induction that for any nonnegative
integer $a$

$$AE^{2}\left(CE^{2}\right)^{a}=AE^{2}=\left[\begin{array}{lll} 1 & 1
& 1 \\
1 & 1
& 1 \\
1 & 1 & 1
\end{array}\right] \phi^{2},$$
therefore the first term in (\ref{equ:trace equation}) is
\begin{align}
\operatorname{tr}\left((AE^{2})^{g-1}(CE^{2})^{-n}\right)&=\operatorname{tr}\left(
\left[\begin{array}{lll} 1 & 1
& 1 \\
1 & 1
& 1 \\
1 & 1 & 1
\end{array}\right]^{g-1} \phi^{2g-2}\right)\nonumber\\
&=3^{g-1}\phi^{2g-2}.
\end{align}
Again induction on nonnegative integers $a,b$ together with simple
calculations imply that
$$(BEC)^{b}(E^{2}C)^{a}=
\left[\begin{array}{lll} 0 & 0
& \frac{t_{1}-t_{2}}{t_{0}-t_{2}} \\
0& 0
& \frac{t_{0}-t_{2}}{t_{1}-t_{0}} \\
0 & 0 & 1
\end{array}\right] \phi^{b},$$
therefore the the second term in (\ref{equ:trace equation}) is
\begin{align}
\operatorname{tr}\left((BEC)^{g-1}(E^{2}C)^{-n}\right)=3^{g-1}\phi^{2g-2}.
\end{align}
Powers of $BCE$ are more difficult to compute, so for computing
the third term in (\ref{equ:trace equation}) we first notice that
$$CEB=\left[\begin{array}{lll} 0 & 0
& 0 \\
0& 0
& 0 \\
3 & 3 & 3
\end{array}\right]\phi^{2},$$
and also $$(ECE)^{a}=ECE$$ for any positive integer $a$, so for
$b>1$, we can write
$$(BCE)^{b}(ECE)^{a}=(B(CEB)^{b-1}CE)(ECE)=3^{b-1}(BCE)(ECE).$$
Easy calculation shows that
$$\operatorname{tr}\left((BCE)(ECE)\right)=\operatorname{tr}(BCE)=
3\phi^{2}.$$ Putting all together, we can find the third term in
(\ref{equ:trace equation}):
\begin{align}\label{equ:the 3rd term}
\operatorname{tr}\left((BCE)^{g-1}(ECE)^{-n}\right)=3^{g-1}\phi^{2g-2}.
\end{align}
By (\ref{equ:trace equation})-(\ref{equ:the 3rd term}), we find
$$Z_{\beta_{cs}}(g\operatorname{|}0,k)=3^{g}\phi^{2g-2},$$
which proves the theorem in this case.

\item[(ii)]  $n\geq 0$\\
We have $U_{2}^{-1}=M+N$, so we can rewrite (\ref{equ:z(g|0,k)})
as
$$Z(g\operatorname{|}0,k)=\operatorname{tr}\left(\left((A+B)(C+E)^{2}
\right)^{g-1}(M+N)^{3n}\right).$$ One can check that
$(M^{2},N)=M^{3}=0$, so $$Z(g\operatorname{|}0,k)=
\operatorname{tr}\left(\left(AE^{2}+B(E,C)+\dots\right)^{g-1}
\left((M,N^{2})+N^{3}\right)^{n}\right).$$ By the same reason as
in the last case we have
\begin{align}\label{equ:z(g|0,k) case2}
&Z_{\beta_{cs}}(g\operatorname{|}0,k)
=\operatorname{tr}\left((AE^{2})^{g-1}(N^{2}M)^{n}\right)
+\operatorname{tr}\left((BEC)^{g-1}(NMN)^{n}\right)\nonumber \\
&\hspace{2.15cm}+\operatorname{tr}\left((BCE)^{g-1}(MN^{2})^{n}\right).
\end{align}
Easy calculations shows that $$BM=MB=0,$$ and for any positive
integer $a$
\begin{align*}
\left(N^{2}M\right)^{a}&=N^{2}M,\\
\left(MN^{2}\right)^{a}&=MN^{2},\\
\left(NMN\right)^{a}&=NMN,\\
(AE^{2})\left(N^{2}M\right)^{a}&=AE^{2},
\end{align*}
and also
\begin{align*}
\operatorname{tr}\left((BCE)(NMN)\right)&=3\phi^{2},\\
\operatorname{tr}\left((BEC)(MN^{2})^{a}\right)&=3\phi^{2}.
\end{align*}
Now calculations similar to the last case together with these
equalities imply that each term in the right hand side of
(\ref{equ:z(g|0,k) case2}) is equal to
$$3^{g-1}\phi^{2g-2},$$ and this proves the theorem in
this case.
\end{enumerate}
For $g=0$, the result is deduced from the semisimplicity of the
TQFT (Corollary \ref{pro:semisimple}, see also Remark
\ref{rem:TQFT}). \qed

\appendix \section{Proof of the Gluing theorem} \label{app:proof of
gluing} We first prove the assertion of Remark \ref{rem:connected
vs disconnected}, which deals with the fact that we don't need to
consider maps with disconnected domains:

\begin{lem} \label{lem:connected vs disconnected}
The contribution of maps with disconnected domain curves in the
section class equivariant Gromov-Witten invariants of the space
$$\mathbb{P}(\mathcal{O}\oplus L_{1}\oplus L_{2})$$ is zero.
\end{lem}
\textsc{Proof:} A disconnected domain curve whose image represents
the class $\beta_{0}+nf$ is a union of some connected components
such that at least the image of one of them represents the class
$n'f$ for a positive $n'$. We have
$$\text{virdim}
\overline{M}(X/\vec{F},n'f)=-(-3H+(2g-2-k_{1}-k_{2})F)\cdot
n'f=3n'>0,$$  so by a discussion similar to Remark \ref{rem:global
theory} one can see that \[
\int_{[\overline{M}_{h}(X/\vec{F},k'f)]^{\text{vir}}} 1=0.
\]
Disconnected invariants can be expressed in terms of the product
of connected ones, and so the lemma follows.\qed

Now we return to the proof of Theorem \ref{subsec:gluing rules}.
We prove the first formula, the proof of the second one is
similar. For simplicity we prove the case $s=0$ and $t=0$.
Extending the argument to the general case is straightforward.

Let $C_{0}$ be a connected curve of genus $g$ with two irreducible
components, $C'$ and $C''$ of genera $g'$ and $g''$ which are
attached together at one point $p$. In other words
$$C_{0}=C'\bigcup_{p=p'=p''} C'',$$ where
$p' \in C'$ and $p'' \in C''$. Now we consider two
$\mathbb{P}^{2}$-bundles
$$X'=\mathbb{P}(\mathcal{O}\oplus L'_{1}\oplus L'_{2})\rightarrow
C',$$ and
$$X''=\mathbb{P}(\mathcal{O}\oplus L''_{1}\oplus L''_{2})\rightarrow
C'',$$ where $L'_{1}$, $L'_{2}$, $L''_{1}$ and $ L''_{2}$ are line
bundles of degrees $k'_{1}$, $k'_{2}$, $k''_{1}$ and $k''_{2}$,
respectively. We attach these two spaces by identifying the
fibers, $F'$ and $F''$ over $p'$ and $p''$, such that the
resulting space is
$$W_{0}=\mathbb{P}(\mathcal{O}\oplus L_{1}\oplus L_{2})\rightarrow
C_{0},$$ where $L_{1}$ and $L_{2}$ are line bundles of degrees
$k_{1}=k'_{1}+k'_{2}$, $k_{2}=k''_{1}+k''_{2}$, respectively. In
other words
$$W_{0}=X'\bigcup_{F=F'=F''}X'',$$ where $F$ is the fiber over
$p$.

Let $W\rightarrow\mathbb{A}^{1}$ be a generic, $1$-parameter
deformation of $W_{0}$ for which the fibers $W_{t}$ for $t\neq0\in
\mathbb{A}^{1}$ are
$$\mathbb{P}(\mathcal{O}\oplus L_{1}\oplus L_{2})\rightarrow C,$$
where $C$ is a smooth curve of genus $g$, and $L_{1}$ and $L_{2}$
are line bundles of degrees $k_{1}$, $k_{2}$.

We follow Jun Li's proof of the degeneration formula in
\cite{Li-relative2}. Let $\W$ be the stack of expanded
degenerations of $W$, with central fiber $\W_{0}$, and let
$\overline{M}_{h}(\W,\beta)$ be the stack of non-degenerate,
pre-deformable, genus $h$, class $\beta$ maps to $\W$, where
$\beta$ is a section class (see \cite{Li-relative1}).

We have the evaluation maps
\[\text{ev}': \overline{M}_{h'}(X'/F',\beta')
\rightarrow F'=F, \] and
\[\text{ev}'': \overline{M}_{h''}(X'/F'',\beta'')
\rightarrow F''=F.\] Li constructs a map
$$\Phi_{\eta}:\overline{M}(X'/F',\beta')\times_{F}
\overline{M}(X''/F'',\beta'')\rightarrow
\overline{M}(\W_{0},\beta),$$ where $\eta$ includes a pair of
classes $(\beta',\beta'')$, such that $\beta=\beta'+\beta''$ and a
pair of genera $(h',h'')$, such that $h=h'+h''$. Then he gives a
virtual cycle formula, which in our case is
\begin{equation} \label{equ:virtual cycle}
[\overline{M}(\W_{0}),\beta]^{vir}=\sum_{\eta}(\Phi_{\eta})_{*}\Delta^{!}
([\overline{M}(X'/F',\beta')]^{vir}\times
[\overline{M}(X''/F'',\beta'')]^{vir}),
\end{equation} where
$\Delta: F\rightarrow F\times F$ is the diagonal map.

\begin{rem}
The torus action on the family $W\rightarrow\mathbb{A}^{1}$ gives
an action on the stack of expanded degeneration, $\W$. One can
check that pre-deformability condition is invariant under this
action, so it induces (canonically) an action on each of the
moduli spaces $\overline{M}_{h'}(X'/F',\beta')$,
$\overline{M}_{h''}(X''/F'',\beta'')$ and
$\overline{M}_{h}(\W,\beta)$. Therefore Li's formula holds in the
equivariant Chow groups.
\end{rem}

If we work with the basis elements $x_{0}$, $x_{1}$, $x_{2}$
(introduced in Section \ref{sec:relative invariants}), for the
equivariant cohomology of the fiber $F$, by using (\ref{equ:basis
relations}), one can see easily that
$$x_{0}^{\vee}=\frac{x_{0}}{T(x_{0})},\hspace{2mm}x_{1}^{\vee}=\frac{x_{1}}{T(x_{1})},
\hspace{2mm}x_{2}^{\vee}=\frac{x_{2}}{T(x_{2})}$$ is its dual
basis, so the cohomology class of the diagonal of $F\times F$ is
given by (see \cite{Mi-St}, Theorem 11.11)
\begin{align*} \label{equ:diagonal}
\operatorname{im}(\Delta)&=\sum_{i=0}^{2}x_{i}\times
x_{i}^{\vee}\nonumber\\
&=x_{0}\times\frac{x_{0}}{T(x_{0})}+x_{1}\times\frac{x_{1}}{T(x_{1})}
+x_{0}\times\frac{x_{2}}{T(x_{2})}.
\end{align*}
Using this, we can rewrite (\ref{equ:virtual cycle}) as
\begin{align*}
&[\overline{M}(\W_{0},\beta)]^{vir}=\\
&\sum_{\eta}(\Phi_{\eta})_{*}\left(\sum_{i=0}^{2}
(x_{i}\cap[\overline{M}(X'/F',\beta')]^{vir}\times
\frac{x_{i}}{T(x_{i})}\cap[\overline{M}(X''/F'',\beta'')]^{vir})\right).
\end{align*}
We now have
\begin{align*}
Z^{h}_{\beta}(g\operatorname{|}k_{1},k_{2})&=\int_{[\overline{M}(\W_{t},\beta)]^{vir}}1\\
&=\int_{[\overline{M}(\W_{0},\beta)]^{vir}}1\\
&=\sum_{\eta}\sum_{i=0}^{2}\int_{[\overline{M}(X'/F',\beta')]^{vir}}x_{i}
\int_{[\overline{M}(X''/F'',\beta'')]^{vir}}\frac{x_{i}}{T(x_{i})}\\
&=\sum_{\eta}\sum_{i=0}^{2}Z^{h'}_{\beta'}(g'\operatorname{|}k'_{1},k'_{2})_{x_{i}}
Z^{h''}_{\beta''}(g''\operatorname{|}k''_{1},k''_{2})^{x_{i}}.
\end{align*}
Then we can write
\begin{align*}
Z(g\operatorname{|}k_{1},k_{2})&=\sum_{\beta\text{ is a section
class}}\sum_{h}u^{2h-2-K_{X}\cdot\beta}Z^{h}_{\beta}(g\operatorname{|}k_{1},k_{2})\\
&=\sum_{\beta}\sum_{h}u^{2h-2-K_{X}\cdot\beta}
\sum_{\eta}\sum_{i}Z^{h'}_{\beta'}(g'\operatorname{|}k'_{1},k'_{2})_{x_{i}}
Z^{h''}_{\beta''}(g''\operatorname{|}k''_{1},k''_{2})^{x_{i}}\\
&=\sum_{\beta,\hspace{1mm}h,\hspace{1mm}\eta,\hspace{1mm}i}
u^{2h'-2-K_{X'}\cdot\beta'}
Z^{h'}_{\beta'}(g'\operatorname{|}k'_{1},k'_{2})_{x_{i}}
u^{2h''-2-K_{X''}\cdot\beta''}
Z^{h''}_{\beta''}(g''\operatorname{|}k''_{1},k''_{2})^{x_{i}}\\
&=\sum_{i}Z(g'\operatorname{|}k'_{1},k'_{2})_{x_{i}}
Z(g''\operatorname{|}k''_{1},k''_{2})^{x_{i}}.
\end{align*}
\section{Some special cases} \label{app:special cases}
In this Appendix, we apply Theorem \ref{thm:matrix formula}, and
drive formulas for some special partition functions of the
equivariant Gromov-Witten invariants of the space
$$\mathbb{P}(\mathcal{O}\oplus L_{1}\oplus L_{2})\rightarrow C,$$
where $C$ is smooth curve of genus $g$.

Similar to Section \ref{sec:proof of theorem}, we use the
following notations for simplicity:
$$
\begin{array} {ll}
A=\left[Z_{\beta_{0}}(1\operatorname{|}0,0)_{x_{a}}^{x_{b}}\right],
&B=\left[Z_{\beta_{0}+f}(1\operatorname{|}0,0)_{x_{a}}^{x_{b}}\right],\\
C_{1}=\left[Z_{\beta_{0}-f}(0\operatorname{|}1,0)_{x_{a}}^{x_{b}}\right],
&C_{2}=\left[Z_{\beta_{0}-f}(0\operatorname{|}0,1)_{x_{a}}^{x_{b}}\right],\\
E_{1}=\left[Z_{\beta_{0}}(0\operatorname{|}1,0)_{x_{a}}^{x_{b}}\right],
&E_{2}=\left[Z_{\beta_{0}}(0\operatorname{|}0,1)_{x_{a}}^{x_{b}}\right],\\
N_{1}=\left[Z_{\beta_{0}}(0\operatorname{|}-1,0)_{x_{a}}^{x_{b}}\right],
&N_{2}=\left[Z_{\beta_{0}}(0\operatorname{|}0,-1)_{x_{a}}^{x_{b}}\right],\\
M_{1}=\left[Z_{\beta_{0}+f}(0\operatorname{|}-1,0)_{x_{a}}^{x_{b}}\right],
&M_{2}=\left[Z_{\beta_{0}+f}(0\operatorname{|}0,-1)_{x_{a}}^{x_{b}}\right].
\end{array}
$$
These matrices are given by (\ref{equ:level creation operators}),
(\ref{equ:level anlhn operators}), (\ref{equ:beta0 g adding}) and
(\ref{equ:beta0+f g adding}).

For proving the theorems in this section, we first write the
formula in Theorem \ref{thm:matrix formula} in each case, in terms
of the matrices above. After expanding it we will get a polynomial
in $t_{0}$, $t_{1}$ and $t_{2}$ (see Remark \ref{rem:global
theory}). We then try to compute the specific term of that
polynomial, which corresponds to the partition function that we
are interested in. All the proofs that we are going to provide is
for the case $g>0$. The case $g=0$ in each theorem, follows from
the semisimplicity of the level $(0,0)$ TQFT (Corollary
\ref{pro:semisimple}, see also Remark \ref{rem:TQFT}).

\begin{thm}
Assume that $k_{1}>0$ and $k_{2}\geq0$. Then the class
$\beta_{0}-k_{1}f$, level $(k_{1},-k_{2})$ equivariant
Gromov-Witten partition function of $X$ is given by
$$
Z_{\beta_{0}-k_{1}f}(g\operatorname{|}k_{1},-k_{2})=(t_{1}-t_{0})^{g+k_{1}-1}
(t_{1}-t_{2})^{g+k_{1}+k_{2}-1}
\left(2\operatorname{sin}\frac{u}{2}\right)^{-2k_{1}-k_{2}}.
$$
\end{thm}
\textsc{Proof:} By Theorem \ref{thm:matrix formula} we have
$$Z(g\operatorname{|}k_{1},-k_{2})=\operatorname{tr}\left((A+B)^{g-1}
(C_{1}+E_{1})^{k_{1}}(M_{2}+N_{2})^{k_{2}}\right).$$ After
expanding this, one can see easily that the only term that
corresponds to the class $\beta_{0}-k_{1}f$ is
$A^{g-1}C_{1}^{k_{1}}N_{2}^{k_{2}}$, or in other words
$$Z_{\beta_{0}-k_{1}f}(g\operatorname{|}k_{1},-k_{2})=
\operatorname{tr}\left(A^{g-1}C_{1}^{k_{1}}N_{2}^{k_{2}}\right).$$
All three matrices $A$, $C_{1}$ and $N_{2}$ are diagonal, so one
can easily compute the right hand side of the above equation:
\begin{align*}
A^{g-1}C_{1}^{k_{1}}N_{2}^{k_{2}}=&\left[\begin{array}{ccc}
(t_{0}-t_{1})(t_{0}-t_{2}) & 0
& 0 \\
0 & (t_{1}-t_{0})(t_{1}-t_{2})
& 0 \\
0 & 0 & (t_{2}-t_{0})(t_{2}-t_{1}))
\end{array}\right]^{g-1}\\
&\left[
\begin{array}{ccc}
0 & 0 & 0 \\
0 &
(t_{1}-t_{0})(t_{1}-t_{2})\phi^{-2} & 0 \\
0 & 0 & 0
\end{array} \right]^{k_{1}}\\
&\left[
\begin{array}{ccc}
(t_{0}-t_{2})\phi^{-1} & 0 & 0 \\
0 &
(t_{1}-t_{2})\phi^{-1} & 0 \\
0 & 0 & 0
\end{array}\right]^{k_{2}}\\
=&\left[\begin{array}{ccc}
0 & 0 & 0 \\
0 &
(t_{1}-t_{0})^{g+k_{1}-1}(t_{1}-t_{2})^{g+k_{1}+k_{2}-1}\phi^{-2k_{1}-k_{2}}
& 0 \\
0 & 0 & 0 \end{array} \right],
\end{align*}
so
$$\operatorname{tr}\left(A^{g-1}C_{1}^{k_{1}}N_{2}^{k_{2}}\right)=
(t_{1}-t_{0})^{g+k_{1}-1}(t_{1}-t_{2})^{g+k_{1}+k_{2}-1}\phi^{-2k_{1}-k_{2}}
.$$ This prove the theorem. \qed

Similar to the theorem above one can prove
\begin{thm}
Assume that $k_{1}\geq0$ and $k_{2}>0$. Then the class
$\beta_{0}-k_{2}f$, level $(-k_{1},k_{2})$ equivariant
Gromov-Witten partition function of $X$ is given by
$$
Z_{\beta_{0}-k_{2}f}(g\operatorname{|}-k_{1},k_{2})=(t_{2}-t_{0})^{g+k_{2}-1}
(t_{2}-t_{1})^{g+k_{1}+k_{2}-1}
\left(2\operatorname{sin}\frac{u}{2}\right)^{-2k_{2}-k_{1}}.
$$
\end{thm}\qed

Also we have the following theorem:
\begin{thm}
Assume that $k>0$. Then the class $\beta_{0}-kf$, level $(k,k)$
equivariant Gromov-Witten partition function of $X$ is given by
\begin{align*}
Z_{\beta_{0}-kf}(g\operatorname{|}k,k)=&\big((t_{1}-t_{0})^{g+k-1}(t_{1}-t_{2})^{g-1}\\
&\hspace{1mm}+(t_{2}-t_{0})^{g+k-1}(t_{2}-t_{1})^{g-1}\big)
\left(2\operatorname{sin}\frac{u}{2}\right)^{-k}.
\end{align*}
\end{thm}
\textsc{Proof:} By Theorem \ref{thm:matrix formula} we have
\begin{align*}
Z(g\operatorname{|}k,k)&=\operatorname{tr}\left((A+B)^{g-1}
(C_{1}+E_{1})^{k}(C_{2}+E_{2})^{k}\right)\\
&=\operatorname{tr}\left((A+B)^{g-1}
(E_{1}E_{2}+C_{1}E_{2}+E_{1}C_{2})^{k}\right),
\end{align*} because $C_{1}C_{2}$=0.

One can now check that
$$Z_{\beta_{0}-kf}(g\operatorname{|}k,k)=\operatorname{tr}\left(A^{g-1}
(C_{1}E_{2}+E_{1}C_{2})^{k}\right).$$ We also have
$$C_{1}E_{2}+E_{1}C_{2}=\left[\begin{array}{ccc}
0 & 0 & 0 \\
0 &
t_{1}-t_{0} & 0 \\
0 & 0 & t_{2}-t_{0}
\end{array} \right]\phi^{-1}.$$ So we can write
\begin{align*}
&\hspace{1.3cm}A^{g-1}(C_{1}E_{2}+E_{1}C_{2})^{k}=\\
&\left[\begin{array}{ccc} (t_{0}-t_{1})(t_{0}-t_{2}) & 0
& 0 \\
0 & (t_{1}-t_{0})(t_{1}-t_{2})
& 0 \\
0 & 0 & (t_{2}-t_{0})(t_{2}-t_{1}))
\end{array}\right]^{g-1}\\
&\left[\begin{array}{ccc}
0 & 0 & 0 \\
0 &
t_{1}-t_{0} & 0 \\
0 & 0 & t_{2}-t_{0}
\end{array} \right]^{k}\phi^{-k}=\\
&\left[\begin{array}{ccc}
0 & 0 & 0 \\
0 &
(t_{1}-t_{0})^{g+k-1}(t_{1}-t_{2})^{g-1} & 0 \\
0 & 0 & (t_{2}-t_{0})^{g+k-1}(t_{2}-t_{1})^{g-1}
\end{array} \right]\phi^{-k}.
\end{align*}
Thus
\begin{align*}
\operatorname{tr}\left(A^{g-1}
(C_{1}E_{2}+E_{1}C_{2})^{k}\right)=&\big((t_{1}-t_{0})^{g+k-1}(t_{1}-t_{2})^{g-1}\\
&\hspace{1mm}+(t_{2}-t_{0})^{g+k-1}(t_{2}-t_{1})^{g-1}\big)
\left(2\operatorname{sin}\frac{u}{2}\right)^{-k}.
\end{align*}
This proves the theorem. \qed

With the same method one can prove the following result:
\begin{thm}
Assume that $k_{1}\geq0$ and $k_{2}\geq0$. Then the class
$\beta_{0}$, level $(-k_{1},-k_{2})$ equivariant Gromov-Witten
partition function of $X$ is given by
\begin{align*}
Z_{\beta_{0}}(g\operatorname{|}&-k_{1},-k_{2})=\\&\begin{cases}
(t_{0}-t_{1})^{g+k_{1}-1}(t_{0}-t_{2})^{g+k_{2}-1}
\left(2\operatorname{sin}\frac{u}{2}\right)^{k_{1}+k_{2}}
&k_{1}>0,\hspace{1mm} k_{2}>0,\\
\big((t_{0}-t_{1})^{g+k_{1}-1}(t_{0}-t_{2})^{g-1}\\
\hspace{1mm}+(t_{2}-t_{0})^{g-1}(t_{2}-t_{1}^{g+k_{1}-1}\big)
\left(2\operatorname{sin}\frac{u}{2}\right)^{k_{1}}
&k_{1}>0,\hspace{1mm} k_{2}=0,\\
\big((t_{0}-t_{1})^{g-1}(t_{0}-t_{2})^{g+k_{2}-1}\\
\hspace{1mm}+(t_{1}-t_{0})^{g-1}(t_{1}-t_{2})^{g+k_{2}-1}\big)
\left(2\operatorname{sin}\frac{u}{2}\right)^{k_{2}}
&k_{1}=0,\hspace{1mm} k_{2}>0,\\
(t_{0}-t_{1})^{g-1}(t_{0}-t_{2})^{g-1}+(t_{1}-t_{0})^{g-1}(t_{1}-t_{2})^{g-1}\\
\hspace{1mm}+(t_{2}-t_{0})^{g-1}(t_{2}-t_{1})^{g-1}
&k_{1}=0,\hspace{1mm} k_{2}=0.
\end{cases}
\end{align*}
\end{thm}\qed

Every theorem that we have already proved in this section is about
computing the partition function for the section class
$\beta_{0}+nf$, where $n$ is the smallest possible integer such
that there is curve in $X$ that represents this class (see Remark
\ref{rem:sum finite} and Lemma \ref{lem:beta=beta+f}). This
partition function is the first term in the sum in Remark
\ref{rem:sum finite}. The next theorem is about computing the
partition function which is the last term in that sum. We restrict
ourself to the level $(0,0)$ case. As we will see, even in this
case the proof is more complicated than the the proof of the
previous theorems in this section.

In the following theorem we are interested in the greatest
possible value of $n$ such that the class $\beta=\beta_{0}+nf$,
level $(0,0)$ partition function is not trivially zero. By Remark
\ref{rem:global theory}, $n$ is the greatest integer that
satisfies
\begin{align*}
d=-K_{X}\cdot \beta&=-(-3H+(2g-2)F)(\beta_{0}+nf)\\
&=3n-(2g-2)<0.
\end{align*}
For $g=3k$, we have $\beta=\beta_{0}+(2k-1)f$ for which
$Z_{\beta}(g\operatorname{|}0,0)$ must be a polynomial of degree
$1$. For $g=3k+1$, we have $\beta=\beta_{0}+2kf$ which is the
Calabi-Yau section class. For $g=3k+2$, we have
$\beta=\beta_{0}+2kf$, and $Z_{\beta}(g\operatorname{|}0,0)$ must
be a polynomial of degree $2$.

\begin{thm}
Let $n$ be the greatest integer that satisfies
$$3n\leq 2g-2,$$ then the class
$\beta_{0}+nf$, level $(0,0)$ equivariant Gromov-Witten partition
function of $X$ is given by
$$Z_{\beta_{0}+nf}(g\operatorname{|}0,0)=\begin{cases} 0 &  g=3k,
\\ 3^{g}\left(2\operatorname{sin}\frac{u}{2}\right)^{2g-2} & g=3k+1,
\\
3^{g-2}(g-1)\big(t_{0}^{2}+t_{1}^{2}+t_{2}^{2}\\-t_{0}t_{1}-t_{0}t_{2}-t_{1}t_{2}\big)
\left(2\operatorname{sin}\frac{u}{2}\right)^{2g-4} & g=3k+2.
\end{cases}
$$
\end{thm}
\textsc{Proof:} Applying Theorem \ref{thm:matrix formula} to this
case, we can write
$$
Z(g\operatorname{|}k,k)=\operatorname{tr}\left((A+B)^{g-1}\right).$$
Now we prove each case in this theorem separately.
\begin{enumerate}
\item[(i)] $g=3k$\\
In this case $n=2k-1$, and one can see that (by using the notation
introduced in Section \ref{sec:proof of theorem})
$$Z_{\beta_{0}+(2k-1)f}(g\operatorname{|}0,0)=\operatorname{tr}\left((A^{k},B^{2k-1})\right).$$
We have $$AB^{2}=\left[\begin{array} {ccc}
1 &1 &1\\
1 &1 &1\\
1 &1 &1\end{array} \right]9\phi^{6},$$ so for any positive integer
$a$ we have
\begin{equation} \label{equ:AB2}
(AB^{2})^{a}=3^{3a-3}\phi^{6a-6}AB^{2}.
\end{equation}
One can prove easily that
\begin{align} \label{equ:B3=0}
B^{3}&=0,\\ \nonumber (ABAB^{2})^{2}&=0.
\end{align}
Now we use the fact that for any two square matrices, $U$ and $V$,
\begin{equation} \label{equ:commutative}
\operatorname{tr}(UV)=\operatorname{tr}(VU).
\end{equation}
Applying this to each term of
$\operatorname{tr}\left((A^{k},B^{2k-1})\right)$ for some number
of times, and using the other equalities above, we can prove that
each of these terms is either zero or is equal to
$$
\operatorname{tr}\left(AB(AB^{2})^{k-1}\right)=3^{3k-6}\phi^{6k-12}
\operatorname{tr}\left(ABAB^{2}\right).$$ However,
$$ABAB^{2}=\left[\begin{array} {ccc}
2t_{0}-t_{1}-t_{2} &2t_{0}-t_{1}-t_{2} &2t_{0}-t_{1}-t_{2}\\
2t_{2}-t_{0}-t_{1} &2t_{1}-t_{0}-t_{2} &2t_{1}-t_{0}-t_{2}\\
2t_{2}-t_{0}-t_{1} &2t_{2}-t_{0}-t_{1}
&2t_{2}-t_{0}-t_{1}\end{array} \right]27\phi^{9},$$ so
$$\operatorname{tr}\left(ABAB^{2}\right)=0.$$ This proves that
each term of $\operatorname{tr}\left((A^{k},B^{2k-1})\right)$ is
zero. So
$$Z_{\beta_{0}+(2k-1)f}(g\operatorname{|}0,0)=0.$$
\item[(ii)] $g=3k+1$\\
In this case $\beta_{0}+nf$ is the Calabi-Yau section class, and
we have proved the theorem for this case in more generality in
Section \ref{sec:proof of theorem}.
\item[(iii)] $g=3k+2$\\
In this case $n=2k$, and this time we have
$$Z_{\beta_{0}+2kf}(g\operatorname{|}0,0)=\operatorname{tr}\left((A^{k+1},B^{2k})\right).$$
The cases $k=0,1$ can be proved by easy calculations, so we assume
that $k>1$. Applying (\ref{equ:commutative}) to each term of this
for some number of times, and using (\ref{equ:B3=0}), we can prove
that each term of $\operatorname{tr}\left((A^{k+1},B^{2k})\right)$
is either zero or is equal to either of $$
\begin{array} {lll}
\operatorname{tr}\left(A(AB^{2})^{k}\right),
&\operatorname{tr}\left((AB)^{2}(AB^{2})^{k-1}\right),
&\operatorname{tr}\left(A^{2}B^{2}(AB^{2})^{k-1}\right),
\end{array}$$
where the number of the terms of the first, the second and the
third kind is $2k+1$, $3k+1$, and  $k$, respectively. So we can
write
\begin{align} \label{equ:Z(g|0,0)}
Z_{\beta_{0}+2kf}(g\operatorname{|}0,0)&=(2k+1)\operatorname{tr}\left(A(AB^{2})^{k}\right)\\
&\hspace{.71mm}+
(3k+1)\operatorname{tr}\left((AB)^{2}(AB^{2})^{k-1}\right)\nonumber\\
&\hspace{.71mm}+\hspace{12.4mm}k\operatorname{tr}\left(A^{2}B^{2}(AB^{2})^{k-1}\right).
\end{align}
One can see that
\begin{align*}
&A(AB^{2})=\\&\left[\begin{array}{ccc} (t_{0}-t_{1})(t_{0}-t_{2})
& (t_{0}-t_{1})(t_{0}-t_{2})
& (t_{0}-t_{1})(t_{0}-t_{2}) \\
(t_{1}-t_{0})(t_{1}-t_{2}) & (t_{1}-t_{0})(t_{1}-t_{2})
& (t_{1}-t_{0})(t_{1}-t_{2}) \\
(t_{2}-t_{0})(t_{2}-t_{1}) & (t_{2}-t_{0})(t_{2}-t_{1}) &
(t_{2}-t_{0})(t_{2}-t_{1}),
\end{array}\right]9\phi^{6},\\
&(AB)^{2}(AB^{2})=\\&\left[\begin{array}{ccc}
1 &1 &1 \\
1 &1 &1 \\
1 &1 &1
\end{array}\right]162(t_{0}^{2}+t_{1}^{2}+
t_{2}^{2}-t_{0}t_{1}-t_{0}t_{2}-t_{1}t_{2})\phi^{12},\\
&A^{2}B^{2}(AB^{2})=\\&\left[\begin{array}{ccc}
(t_{0}-t_{1})(t_{0}-t_{2}) & (t_{0}-t_{1})(t_{0}-t_{2})
& (t_{0}-t_{1})(t_{0}-t_{2}) \\
(t_{1}-t_{0})(t_{1}-t_{2}) & (t_{1}-t_{0})(t_{1}-t_{2})
& (t_{1}-t_{0})(t_{1}-t_{2}) \\
(t_{2}-t_{0})(t_{2}-t_{1}) & (t_{2}-t_{0})(t_{2}-t_{1}) &
(t_{2}-t_{0})(t_{2}-t_{1}),
\end{array}\right]243\phi^{12}.
\end{align*}
Using (\ref{equ:AB2}), we can write
\begin{align*}
\operatorname{tr}\left(A(AB^{2})^{k}\right)&=3^{3k-3}\phi^{6k-6}
\operatorname{tr}\left(A(AB^{2})\right)\\
&=3^{3k-1}(t_{0}^{2}+t_{1}^{2}+t_{2}^{2}-t_{0}t_{1}-t_{0}t_{2}-t_{1}t_{2})\phi^{6k},
\end{align*}
\begin{align*}
\operatorname{tr}\left((AB)^{2}(AB^{2})^{k-1}\right)&=3^{3k-6}\phi^{6k-12}
\operatorname{tr}\left((AB)^{2}(AB^{2})^{k-1}\right)\\
&=2\cdot3^{3k-1}(t_{0}^{2}+t_{1}^{2}+t_{2}^{2}-t_{0}t_{1}-t_{0}t_{2}-t_{1}t_{2})\phi^{6k},
\end{align*}
\begin{align*}
\operatorname{tr}\left(A^{2}B^{2}(AB^{2})^{k-1}\right)&=3^{3k-6}\phi^{6k-12}
\operatorname{tr}\left(A^{2}B^{2}(AB^{2})^{k-1}\right)\\
&=3^{3k-1}(t_{0}^{2}+t_{1}^{2}+t_{2}^{2}-t_{0}t_{1}-t_{0}t_{2}-t_{1}t_{2})\phi^{6k},
\end{align*}
We define
$$Q:=t_{0}^{2}+t_{1}^{2}+t_{2}^{2}-t_{0}t_{1}-t_{0}t_{2}-t_{1}t_{2}.$$
Putting these equalities into (\ref{equ:Z(g|0,0)}), we get
\begin{align*}
Z_{\beta}(g\operatorname{|}0,0)&=3^{3k-1}(2k+1+2(3k+1)+k)Q\phi^{6k}\\
&=3^{3k}(3k+1)Q\phi^{6k}\\
&=3^{g-2}(g-1)Q\phi^{2g-4},
\end{align*} and this proves the theorem in this case.
\end{enumerate}\qed

\end{document}